\documentclass [a4paper, 10pt] {article}
\usepackage {amsmath}
\usepackage {amssymb}
\usepackage {latexsym}
\usepackage[dvips]{graphicx}
\usepackage {fancyhdr}
\newtheorem {theorem}{Theorem}[section]
\newtheorem {Corollary}[theorem]{Corollary}
\newtheorem {Definition}[theorem]{Definition}
\newtheorem {Lemma}[theorem]{Lemma}
\newtheorem {Remark}[theorem]{Remark}
\newtheorem {Proposition}[theorem]{Proposition}

\begin{document}

\title{\LARGE $L^{\infty}$ estimates and integrability by compensation in Besov-Morrey spaces and applications} 

\author{Laura Gioia Andrea Keller \footnote{Department of Mathematics, ETH Z\"urich, 8092 Z\"urich, Switzerland}}
\maketitle

\begin{abstract}
$L^{\infty}$ estimates in the integrability by compensation result of H. Wente (\cite{W}) fail in dimension larger than two when Sobolev spaces are replaced by the ad-hoc Morrey spaces (in dimension $n\geq 3$). \\
However, in this paper we prove that $L^{\infty}$ estimates hold in arbitrary dimension when Morrey spaces are replaced by their Littlewood Paley counterparts: Besov-Morrey spaces. \\
As an application we prove the existence of conservation laws to solution of elliptic systems of the form 
\begin{displaymath}
-\Delta u= \Omega \cdot \nabla u
\end{displaymath}
where $\Omega$ is antisymmetric and both $\nabla u$ and  $\Omega$ belong to these Besov-Morrey spaces for which the system is critical.
\end{abstract}


\section{Introduction}

In this section we will give the precise statement of our results and add some remarks. \\
For the sake of simplicity, in what follows we will use the abbreviation $a_x$ for $\frac{\partial}{\partial x}a$.\\

Our work was motivated by Rivi\`ere's \cite{Riv} article about Schr\"odinger systems with antisymmetric potentials, i.e. systems of the form
\begin{equation}\label{Schroedinger}
-\Delta u = \Omega \cdot \nabla u
\end{equation}
with $u \in W^{1,2}(\omega, \mathbb{R}^m)$ and $\Omega \in L^2(\omega, so(m)\otimes \Lambda^1\mathbb{R}^n)$, $\omega \subset \mathbb{R}^n$.\\
The differential equation (\ref{Schroedinger}) has to be understood in the following sense: \\
For all indices $i \in \left \{ 1, \dots , m \right \}$ we have $- \Delta u^i = \sum_{j=1}^m\Omega^i_j \cdot \nabla u^j$ and $L^2(\omega, so(m) \otimes \Lambda^1\mathbb{R}^n)$ means that $\forall \; i,j \in \left \{ 1, \dots , m \right \}, \; \Omega^i_j \in L^2(\omega,\Lambda^1\mathbb{R}^n) \; \textrm{and} \; \Omega^i_j=-\Omega^j_i$.\\
In particular, it was the result that in dimension $n=2$ solutions to (\ref{Schroedinger}) are continuous which attracted our interest.\\

The interest for such systems originates in the fact that they "encode'' all Euler-Lagrange equations for conformally invariant quadratic Lagrangians in dimension 2 (see \cite{Riv} and also \cite{Gr}).

In what follows we will take $\omega=B^n_1(0)$, the n-dimensional unit ball.\\

In the above cited work, there were three crucial ideas: 
\begin{itemize}
\item
\textbf{Antisymmetry of $\Omega$} \\
If we drop the assumption that $\Omega$ is symmetric, there may occur solutions which are not continuous as the following example shows:\\
Let $n=2$, $u^i= 2 \log \log \frac{1}{r}$ for $i=1,2$ and let 
\begin{displaymath}
\Omega = \left ( \begin{array}{cc}
\nabla u^1 & 0  \\
0 & \nabla u^2  
\end{array} \right )
\end{displaymath}
Obviously, $u$ satisfies equation (\ref{Schroedinger}) with the given $\Omega$ but is not continuous.
\item
\textbf{Construction of conservation laws}\\
In fact, once there exist $A \in L^{\infty}(B^n_1(0),M_m(\mathbb{R}))\cap W^{1,2}(B^n_1(0),M_m(\mathbb{R}))$ such that
\begin{equation}\label{GlA}
d^{\ast}(dA-A\Omega)=0.
\end{equation}
for given $\Omega \in L^2(B^n_1(0), so(m) \otimes \Lambda^1\mathbb{R}^n)$, then any solution $u$ of (\ref{Schroedinger}) satisfies the following conservation law
\begin{equation}\label{conslaw}
d(\ast A du + (-1)^{n-1}(\ast B)\wedge du)  =0
\end{equation}
where $B$ satisfies $-d^{\ast}B=dA-A\Omega$.\\
The existence of such an $A$ (and $B$) is proved by Rivi\`ere in \cite{Riv} and relies on a \textbf{non linear Hodge decomposition} which can also be interpreted as a \textbf{change of gauge}. (see in our case theorem \ref{ndimRiv})
\item
\textbf{Understanding the linear problem}\\
The proof of the above mentioned regularity result uses the result below for the linear problem:
\begin{theorem}
(\cite{W},\cite{CLMS}, \cite{Tar})\\
Let $a,b$ satisfy $\nabla a, \nabla b \in L^2$ and let $\varphi$ be the unique solution to
\begin{equation} \label{Dirichlet}
\left \{ \begin{array}{ll}
-\Delta \varphi = \nabla a \cdot \nabla^{\perp} b = \ast (da \wedge db) =a_x b_y - a_y b_x \; \textrm{in $B^n_1(0)$} \\
\varphi =0 \; \textrm{on $\partial B^n_1(0)$}. 
\end{array} \right .
\end{equation}
Then $\varphi$ is continuous and it holds that
\begin{equation}\label{Absch}
\vert \vert \varphi \vert \vert_{\infty} + \vert \vert \nabla \varphi \vert \vert_{2} + \vert \vert \nabla^2 \varphi \vert \vert_1 
\leq C \vert \vert \nabla a \vert \vert _2 \; \vert \vert \nabla b \vert \vert_2.
\end{equation}
Note that the $L^{\infty}$ estimate in (\ref{Absch}) is the key point for the existence of $A$, $B$ satisfying (\ref{GlA}).
\end{theorem}

\end{itemize}
A more detailed explanation of these key points and their interplay can be found in Rivi\`ere's overview \cite{RivLect}.\\
 
Our strategy to extend the cited regularity result to domains of arbitrary dimension is to find first of all a good generalisation of Wente's estimate. Here, the first question is to detect a suitable substitute for $L^2$ since obviously for $n\geq 3$ from the fact that $a,b \in W^{1,2}$ we can not conclude that $\varphi$ is continuous.  So we have to reduce our interest to a smaller space than $L^2$. A first idea is to look at the Morrey space $\mathcal{M}^n_2$, i.e. at the spaces of all functions $f \in L^2_{loc}(\mathbb{R}^n)$ such that
\begin{displaymath}
\vert \vert f \vert \mathcal{M}^n_2 \vert \vert = \sup_{x_0 \in \mathbb{R}^n } \sup_{R>0} R^{1-n/2} \vert \vert f \vert L^2(B(x_0,R))\vert \vert < \infty.
\end{displaymath}
The choice of this space was motivated by the following observation (for details see \cite{RivSt}):\\
For stationary harmonic maps $u$ we have the following monotonicity estimate
\begin{displaymath}
r^{2-n}\int_{B^n_r(x_0)}\vert \nabla u \vert^2  \leq R^{2-n}\int_{B^n_R(x_0)}\vert \nabla u \vert^2
\end{displaymath}
for all $r\leq R$. From this, it is rather natural to look at the Morrey space $\mathcal{M}^n_2$.\\

Unfortunately, this first try is not successful as the following counterexample in dimension $n=3$ shows:\\
Let $a=\frac{x_1}{\vert x \vert}$ and $b=\frac{x_2}{\vert x \vert}$. As required $\nabla a, \nabla b \in \mathcal{M}^3_2(B^3_1(0))$. The results in (\cite{CLMS}) imply that the unique solution $\varphi$ of (\ref{Dirichlet}) satisfies $\nabla ^2 \varphi \in \mathcal{M}^{\frac{3}{2}}_1$  but $\varphi$ is not bounded!\\
Therefore, in \cite{RS} the attempt to construct conservation laws for (\ref{Schroedinger}) in the framework of Morrey spaces fails.\\

Another drawback is that $C^{\infty}$ is not dense in $\mathcal{M}^n_2$.
This point is particularly important if one has in mind the proof via paraproducts of Wente's $L^{\infty}$ bound for the solution $\varphi$.\\
In this paper we shall study $L^{\infty}$ estimates by replacing the Morrey spaces $\mathcal{M}^n_2$ by their "nearest'' Littlewood Paley counterpart, the  Besov-Morrey spaces $B^0_{\mathcal{M}^n_2,2}$, i.e. the spaces of $f \in \mathcal{S}'$ such that
\begin{displaymath}
\Big(\sum_{j=0}^{\infty} \vert \vert \mathcal{F}^{-1}\varphi _j \mathcal{F}f \vert \mathcal{M}^n_2(\mathbb{R}^n)\vert \vert^2\Big)^{\frac{1}{2}} < \infty
\end{displaymath}
where $\varphi=\left \{ \varphi_j \right \}_{j=0}^{\infty}$ is a suitable partition of unity.\\
It turns out that we have a suitable density result at hand, see lemma \ref{density}.
These spaces were introduced by Kozono and Yamazaki in \cite{KY} and applied to the study of the Cauchy problem for the Navier-Stokes equation and semilinear heat equation (see also \cite{Maz2}).\\

Note, that we have the following \textbf{natural embeddings}: $B^0_{\mathcal{M}^n_2,2} \subset \mathcal{M}^n_2$ (see lemma \ref{inMorrey}) and on compact subsets $B^0_{\mathcal{M}^n_2,2}$ is a natural subset of $L^2$ (see lemma \ref{aufKomp}).\\

The success to which these Besov-Morrey spaces give rise relies crucially on the fact that \textbf{we first integrate and then sum}! \\
In the spirit of the scales of Triebel-Lizorkin and Besov spaces (definition are restated in the next section) where we have for $0<q\leq \infty$ and $0<p<\infty$
\begin{displaymath}
B^s_{p,\min\left \{ p,q \right \}} \subset F^s_{p,q} \subset B^s_{\max \left \{ p,q \right \}}
\end{displaymath}
and due to the fact that for $1<q\leq p <\infty$
\begin{displaymath}
\vert \vert f \vert \vert _{\mathcal{M}^p_q} \simeq \Big \vert \Big \vert \Big (\sum_{j=0}^{\infty} \vert \mathcal{F}^{-1}\varphi _j \mathcal{F}f \vert^2  \Big )^{\frac{1}{2}} \Big \vert \Big \vert_{\mathcal{M}^p_q}
\end{displaymath}
it is obvious to exchange the order of summability and integrability in order to find a smaller space starting from a given one.\\

A more detailed exposition of the framework of Besov-Morrey spaces is given in the next section.\\

We have

\begin{theorem}\label{regndim}
\begin{itemize}
\item[i)]
Assume that $a,b \in B^0_{\mathcal{M}^n_2,2}$, and assume further that
\begin{displaymath}
a_x, a_y,b_x,b_y \in B^0_{\mathcal{M}^n_2,2}\; \textrm{where} \; x,y=z_i,z_j \: \textrm{with} \; i,j, \in \left \{ 1, \dots , n \right \} .
\end{displaymath}
Then any solution of
\begin{equation}
- \Delta u= a_xb_y-a_yb_x \nonumber
\end{equation}
is continuous and bounded.
\item[ii)]
Assume that $a_x,a_y,b_x$ and $b_y$ are distributions whose support is contained in $B^n_1(0)$ and belong to $B^0_{\mathcal{M}^n_2,2}$, $n\geq3$.\\ Moreover, let $u$ be a solution ( in the sense of distributions ) of
\begin{displaymath}
-\Delta u = a_xb_y-b_xa_y.
\end{displaymath}
Then it holds
\begin{displaymath}
\nabla u \in B^0_{\mathcal{M}^n_2,1}.
\end{displaymath}
\item[iii)]
Assume that $a_x,a_y,b_x$ and $b_y$ are distributions whose support in $B^n_1(0)$ and belong to $B^0_{\mathcal{M}^n_2,2}$. \\
Moreover, let $u$ be a solution ( in the sense of distributions ) of
\begin{displaymath}
-\Delta u = a_xb_y-b_xa_y.
\end{displaymath}
Then it holds
\begin{displaymath}
\nabla^2 u \in  B^{-1}_{\mathcal{M}^n_2,1} \subset B^{-2}_{\infty,1}.
\end{displaymath}
\end{itemize}
\end{theorem}

\begin{Remark}
\textnormal{
\begin{itemize}
\item
If we reduce our interest to dimension $n=2$, our assumption in the theorem below coincide with the original ones in Wente's framework due to the fact that $\mathcal{M}^2_2=L^2$ and $B^0_{2,2}=L^2=F^0_{2,2}$.
\item
Obviously we have the a-priori bound
\begin{displaymath}
\vert \vert u \vert \vert_{\infty} \leq C \big ( \vert \vert a \vert B^0_{\mathcal{M}^n_2,2} \vert \vert + \vert \vert \nabla a \vert B^0_{\mathcal{M}^n_2,2} \vert \vert \big) \big (\vert \vert b \vert B^0_{\mathcal{M}^n_2,2} \vert \vert + \vert \vert  \nabla b  \vert B^0_{\mathcal{M}^n_2,2} \vert \vert \big).
\end{displaymath}
\item
Now, if we use a homogeneous partition of unity instead of an inhomogeneous as before, our result holds if we replace the spaces $B^0_{\mathcal{M}^n_2,2}$ by the spaces $\mathcal{N}^0_{n,2,2}$. For further information about these homogeneous function spaces we refer to Mazzucato's article \cite{Maz2}.
\item
Note that the estimate $\nabla u \in B^0_{\mathcal{M}^n_2,1}$ implies that $u$ is bounded and continuous.
\end{itemize}
}
\end{Remark}

As an application of what we did so far, we would like to present an adaptation of Rivi\`ere's construction of conservation laws via gauge transformation (see \cite{Riv}) to our setting, more precisely we are able to prove the following assertion:

\begin{theorem}\label{gauge}
Let $n\geq 3$. There exist constants $\varepsilon(m)>0$ and $C(m) >0$ such that for every $\Omega \in B^0_{\mathcal{M}^n_2,2}(B^n_1(0), so(m)\otimes \Lambda^1\mathbb{R}^n)$ which satisfies
\begin{displaymath} 
\vert \vert \Omega \vert B^0_{\mathcal{M}^n_2,2} \vert \vert \leq \varepsilon(m)
\end{displaymath}
there exist $A\in L^{\infty}(B^n_1(0),Gl_m(\mathbb{R}))\cap B^1_{\mathcal{M}^n_2,2}$ and $B \in B^1_{\mathcal{M}^n_2,2}(B^n_1(0),M_m(\mathbb{R})\otimes \Lambda^2\mathbb{R}^n)$ such that
\begin{itemize}
\item[i)]
\begin{displaymath}
d_{\Omega}:= dA-A\Omega = -d^{*}B=-*d*B
\end{displaymath}
\item[ii)]
\begin{displaymath}
\vert \vert \nabla A \vert B^0_{\mathcal{M}^n_2,2} \vert \vert +\vert \vert \nabla A^{-1} \vert B^0_{\mathcal{M}^n_2,2} \vert \vert +
\int_{B^n_1(0)} \vert \vert \textrm{dist}(A,SO(m))\vert \vert_{\infty}^2 \leq C(M) \vert \vert \Omega \vert B^0_{\mathcal{M}^n_2,2} \vert \vert
\end{displaymath}
\item[iii)]
\begin{displaymath}
\vert \vert \nabla B \vert B^0_{\mathcal{M}^n_2,2} \vert \vert \leq C(m)\vert \vert \Omega \vert B^0_{\mathcal{M}^n_2,2} \vert \vert.
\end{displaymath}
\end{itemize}
\end{theorem}

This finally leads to the following regularity result:

\begin{Corollary}\label{ndimRiv}
Let the dimension $n$ satisfy $n\geq 3$. Let $\varepsilon(m)$, $\Omega$, $A$ and $B$ be as in theorem \ref{gauge}.
Then any solution $u$ of 
\begin{displaymath}
- \Delta u = \Omega \cdot \nabla u
\end{displaymath}
satisfies the conservation law 
\begin{displaymath}
d(\ast A du + (-1)^{n-1}(\ast B)\wedge du)  =0.
\end{displaymath}
Moreover, any distributional solution of $\Delta u = -\Omega \cdot \nabla u$
which satisfies in addition
\begin{displaymath}
\nabla u \in B^0_{\mathcal{M}^n_2,2}
\end{displaymath}
is continuous.
\end{Corollary}

\begin{Remark}
\textnormal{
Note that the continuity assertion of the above corollary is already contained in \cite{RivSt}, but our result differs from \cite{RivSt} (see also \cite{Schi} for a modification of the proof of Rivi\`ere and Struwe)  in so far, as on one hand we do not impose any smallness of the norm of the gradient of a solution and really construct $A$ and $B$ (see theorem \ref{gauge}) and not only construct $\Omega$ and $\xi$ such that $P^{-1} d P+P^{-1}\Omega P=*d\xi$, but on the other hand work in a slightly smaller space.
}
\end{Remark}

The present article is organised as follows: After recalling some basic definitions and preliminary facts in section 2 we give in the third section the proofs of the statements claimed before. \\

\textbf{Acknowledgement}
The author would like to thank Professor T. Rivi\`ere for having pointed out the present problem to her and for his support.

\section{Definitions and preliminary results}

We recall the important definitions and state basic results we will use.

\subsection{Besov and Triebel-Lizorkin spaces}

\subsubsection{Non-homogeneous Besov and Triebel-Lizorkin spaces}

In order to define them we have to introduce some additional notions:

\begin{Definition}[$\mathbf{\Phi(\mathbb{R}^n)}$]
Let $\mathbf{\Phi (\mathbb{R}^n)}$ be the collection of all systems \\
$\varphi = \left \{ \varphi _j(x) \right\}_{j=0}^{\infty} \subset \mathcal{S}(\mathbb{R}^n)$ such that 
\begin {displaymath}
\left\{ \begin{array}{ll}
supp \; \varphi _0 \subset \left\{x \vert \; \vert x \vert \leq 2\right\} & \\
supp \; \varphi _j \subset \left\{x \vert \; 2^{j-1}\leq \vert x \vert \leq 2^{j+1} \right\} & \textrm{if $j=1,2,3,\dots , $}
\end{array} \right.
\end {displaymath}
for every multi-index $\alpha $ there exists a positive number $C_{\alpha }$ such that
\begin{displaymath}
2^{j\vert\alpha \vert} \vert D^{\alpha }\varphi _j (x) \vert \leq C_{\alpha} \; \textrm{for all $j=1,2,3,\dots \;$ and all
$x \in \mathbb{R}^n$}
\end{displaymath}
and
\begin{displaymath}
\sum_{j=0}^{\infty} \varphi _j(x) =1 \; \forall x \in \mathbb{R}^n
\end{displaymath}
\end{Definition}

\begin{Remark}\label{example}
\textnormal{
\begin{itemize}
\item
Note that in the above expression $\sum_{j=0}^{\infty} \varphi _j(x) =1$ the sum is locally finite!
\item
Example of a system $\varphi$ which belongs to $\Phi (\mathbb{R}^n)$:\\
We start with an arbitrary $C^{\infty}_0(\mathbb{R}^n)$ function $\psi  $ which has the following properties: 
$\psi (x ) =1$ for $ \vert x \vert \leq 1$ and $\psi (x ) =0$ for $\vert x \vert \geq \frac{3}{2}$. We set $\varphi _0(x )=\psi(x )$, 
$\varphi _1(x )=\psi(\frac{x}{2})-\psi(x) $,   and 
$\varphi _j(x )=\varphi_1 (2^{-j+1}x )$, $j\geq 2$. Then it is easy to check that this family $\varphi $ 
satisfies the requirements of our definition.  \\
Moreover, we have $\sum_{j=0}^{n}\varphi_j(x)=\psi(2^{-n}x)$, $n\geq 0$.\\
By the way, other examples of $\varphi \in \Phi $,
apart from this one, can be found in \cite{RS}, \cite {T} or \cite{C}. ) 
\end{itemize}
}
\end{Remark}

Now, we can state the definitions of the above mentioned Besov and Triebel-Lizorkin spaces.

\begin{Definition}[Besov spaces and Triebel-Lizorkin spaces]
Let $ -\infty < s < \infty$, let $0< q \leq\infty$ and let $\varphi \in \Phi (\mathbb{R}^n)$.
\begin{itemize}
\item[i)]
If $0<p\leq \infty$ then
the (non-homogeneous) Besov spaces $\mathbf{B^s_{p,q}(\mathbb{R}^n)}$ consist of all $f \in \mathcal{S}'$ such that the following inequality holds
\begin{displaymath}
\vert \vert f \vert B^s_{p,q}(\mathbb{R}^n)\vert \vert ^{\varphi }= 
\vert \vert 2^{js}\mathcal{F}^{-1}\varphi _j \mathcal{F}f \vert l^q(L^p(\mathbb{R}^n))\vert\vert < \infty
\end{displaymath}
\item[ii)]
If $0<p<\infty$ then
the (non-homogeneous) Triebel-Lizorkin spaces $\mathbf{F^s_{p,q}(\mathbb{R}^n)}$ consist of all $f \in \mathcal{S}'$ such that the following inequality holds
\begin{displaymath}
\vert \vert f \vert F^s_{p,q}(\mathbb{R}^n)\vert \vert ^{\varphi }= 
\vert \vert 2^{js}\mathcal{F}^{-1}\varphi _j \mathcal{F}f \vert L^p(\mathbb{R}^n,l^q)\vert\vert < \infty
\end{displaymath} 
\item[iii)]
If $p=\infty$ then
the spaces $\mathbf{F^s_{\infty,q}(\mathbb{R}^n)}$ consist of all $f \in \mathcal{S}'$ such that \\
$\exists \left\{f_k(x)\right\} _{k=0}^{\infty} \subset L^{\infty}(\mathbb{R}^n)$ such that the following holds
\begin{displaymath}
f=\sum_{k=0}^{\infty} \mathcal{F}^{-1}\varphi _k\mathcal{F}f_k \; \textrm{in} \; \mathcal{S}'(\mathbb{R}^n)
\end{displaymath}
and 
\begin{displaymath}
\vert \vert 2^{sk}f_k  \vert  L^{\infty}(\mathbb{R}^n,l^q) \vert \vert < \infty.
\end{displaymath}
Moreover we set
\begin{displaymath}
\vert \vert f \vert F^s_{\infty,q}(\mathbb{R}^n)\vert \vert ^{\varphi }= 
\inf \vert \vert 2^{sk} f_k  \vert L^{\infty}(\mathbb{R}^n,l^q)\vert\vert
\end{displaymath} 
where the infimum is taken over all admissible representations of $f$.
\end{itemize}
Here $\mathcal{F}$ denotes the Fourier transform and
\begin{displaymath}
\vert \vert f_k \vert l^q(L^p(\mathbb{R}^n))\vert \vert = 
\Bigg( \sum_{k=0}^{\infty} \big( \int \vert f_k(x)\vert^p dx \big)^{\frac{q}{p}}\Bigg)^{\frac{1}{q}}
\end{displaymath}
and
\begin{displaymath} 
\vert \vert f_k \vert L^p(\mathbb{R}^n,l^q)\vert \vert =  
\Bigg( \int \big( \sum_{k=0}^{\infty} \vert f_k(x)\vert^q \big)^{\frac{p}{q}}dx\Bigg)^{\frac{1}{p}}.
\end{displaymath}
\end {Definition}

Recall that the spaces $B^s_{p,q}$ and $F^s_{p,q}$ are independent of the choice of $\varphi $ (see \cite{T}).\\

Most of the important fact (embeddings, relation with other function spaces, multiplier assertions and so on) about these spaces can be found in \cite{RS} and \cite{T}. In what follows we will give precise indications where a result we use is proved.

\subsubsection{Besov-Morrey spaces}

In stead of combining $L^p$-norms ans $l^q$-norm one can also combine Morrey- (respectively Morrey-Campanato-) norms with $l^q$-norms. This idea was first introduced and applied by  Kozono and Yamazaki in \cite{KY}.\\
In order to make the whole notation clear and to avoid misunderstanding we will recall some definitions. \\
We start with the definition of Morrey spaces

\begin{Definition}[Morrey spaces]
Let $1 \leq q \leq p < \infty$.
\begin{itemize}
\item[i)]
The Morrey spaces $\mathbf{\mathcal{M}_q^p(\mathbb{R}^n)}$ consist of all $f \in L^q_{loc}(\mathbb{R}^n)$ such that
\begin{displaymath}
\vert \vert f \vert \mathcal{M}^p_q \vert \vert = \sup_{x_0 \in \mathbb{R}^n } \sup_{R>0} R^{n/p-n/q} \vert \vert f \vert L^q(B(x_0,R))\vert \vert < \infty 
\end{displaymath}
\item[ii)]
The local Morrey spaces $\mathbf{M_q^p(\mathbb{R}^n)}$ consist of all $f \in L^q_{loc}(\mathbb{R}^n)$ such that
\begin{displaymath}
\vert \vert f \vert \mathcal{M}^p_q \vert \vert = \sup_{x_0 \in \mathbb{R}^n } \sup_{0<R\leq 1} R^{n/p-n/q} \vert \vert f \vert L^q(B(x_0,R))\vert \vert < \infty
\end{displaymath}
where $B(x_0,R)$ denotes the closed ball in $\mathbb{R}^n$ with center $x_0$ and radius $R$.
\end{itemize}
\end{Definition}

Note that it is easy to see that the spaces $\mathcal{M}^p_q$ and $M^p_q$ coincide on compactly supported functions.

Apart from these spaces of regular distributions, i.e. function belonging to $L^1_{loc}$, in the case $q=1$ we are even allowed to look at measures in stead of functions. More precisely we have the following measure spaces of Morrey type. They will become useful later on in a rather technical context.

\begin{Definition}[Measure spaces of Morrey type]
Let $1 \leq p < \infty$.
\begin{itemize}
\item[i)]
The measure spaces of Morrey type $\mathbf{\mathcal{M}^p(\mathbb{R}^n)}=\mathcal{M}^p$ consist of all Radon measures $\mu$ such that
\begin{displaymath}
\vert \vert \mu \vert \mathcal{M}^p \vert \vert = \sup_{x_0 \in \mathbb{R}^n } \sup_{R>0} R^{n/p-n} \vert \mu \vert (B(x_0,R))< \infty.
\end{displaymath}
\item[ii)]
The local measure spaces of Morrey type $\mathbf{M_q(\mathbb{R}^n)}=M^p$ consist of all Radon measures $\mu$ such that
\begin{displaymath}
\vert \mu \vert \mathcal{M}^p \vert \vert = \sup_{x_0 \in \mathbb{R}^n } \sup_{0<R\leq 1} R^{n/p-n} \vert \mu \vert(B(x_0,R)) < \infty
\end{displaymath}
where as above $B(x_0,R)$ denotes the closed ball in $\mathbb{R}^n$ with center $x_0$ and radius $R$.
\end{itemize}
\end{Definition}

Remember that all the spaces we have seen so far, i. e. $\mathcal{M}^p_q$, $M^p_q$, $\mathcal{M}^p$ and $M^p$ are Banach spaces with the norms indicated before. Moreover, $\mathcal{M}^p_1$ and $M^p_1$ can be considered as closed subspaces of $\mathcal{M}^p$ and $M^p$ respectively, consisting of all those measures which are absolutely continuous with respect to the Lebesgue measure. \\
For details, see e.g. \cite{KY}.

Once we have the above definition of Morrey spaces (of regular distributions), we now define the Besov-Morrey spaces in the same way as we constructed the Besov spaces, of course with the necessary changes.

\begin{Definition}[Besov-Morrey spaces]
Let $1\leq q \leq p < \infty$, $1\leq r \leq \infty$ and $s \in \mathbb{R}$.
\begin{itemize}
\item[i)]
Let $\varphi \in \dot{\Phi} (\mathbb{R}^n)$. The homogeneous Besov-Morrey spaces $\mathbf{\mathcal{N}^s_{p,q,r}}$ consist of all $f \in \mathcal{Z}'$ such that
\begin{displaymath}
\vert \vert f \vert \mathcal{N}^s_{p,q,r}(\mathbb{R}^n)\vert \vert ^{\varphi }= 
\Big(\sum_{j= - \infty}^{\infty}2^{jsr} \vert \vert \mathcal{F}^{-1}\varphi _j \mathcal{F}f \vert \mathcal{M}^p_q(\mathbb{R}^n)\vert \vert^r\Big)^{\frac{1}{r}} < \infty.
\end{displaymath}
\item[ii)]
Let $\varphi \in \Phi (\mathbb{R}^n)$. The inhomogeneous Besov-Morrey spaces $\mathbf{N^s_{p,q,r}}$ consist of all $f \in \mathcal{S}'$ such that
\begin{displaymath}
\vert \vert f \vert N^s_{p,q,r}(\mathbb{R}^n)\vert \vert ^{\varphi }= 
\Big(\sum_{j= 0}^{\infty}2^{jsr} \vert \vert \mathcal{F}^{-1}\varphi _j \mathcal{F}f \vert M^p_q(\mathbb{R}^n)\vert \vert^r\Big)^{\frac{1}{r}} < \infty.
\end{displaymath}
\end{itemize}
\end{Definition}

Note  that since $L^p(\mathbb{R}^n)=\mathcal{M}^p_p(\mathbb{R}^n)$ the framework of the $\mathcal{N}^s_{p,q,r}(\mathbb{R}^n)$ can be seen as a generalisation of the framework of the homogeneous Besov spaces. 

In our further work we will crucially use still another variant of spaces which are defined via Paley-Littlewood decomposition. We will use the decomposition into frequencies of positive power but measure the single contributions in a homogeneous Morrey norm:

\begin{Definition}[The spaces $B^s_{\mathcal{M}^p_q,r}$]
\begin{itemize}
\item[i)]
Let $1\leq q \leq p < \infty$, $1\leq r \leq \infty$ and $s \in \mathbb{R}$.
Let $\varphi \in \Phi (\mathbb{R}^n)$. The spaces $\mathbf{B^s_{\mathcal{M}^p_q,r}}$ consist of all $f \in \mathcal{S}'$ such that
\begin{displaymath}
\vert \vert f \vert B^s_{\mathcal{M}^p_q,r}(\mathbb{R}^n)\vert \vert ^{\varphi }= 
\Big(\sum_{j=0}^{\infty}2^{jsr} \vert \vert \mathcal{F}^{-1}\varphi _j \mathcal{F}f \vert \mathcal{M}^p_q(\mathbb{R}^n)\vert \vert^r\Big)^{\frac{1}{r}} < \infty.
\end{displaymath}
\item[ii)]
The spaces $\mathbf{B^s_{\mathcal{M}^p_q,r}}(\Omega)$ where $\Omega$ is a bounded domain in $\mathbb{R}^n$ consist of all $f \in B^s_{\mathcal{M}^p_q,r}$ which in addition have compact support contained in $\Omega$.
\end{itemize}
\end{Definition}

\begin{Remark}
\textnormal{
\begin{itemize}
\item[i)]
Again, as in the case of Besov and Triebel-Lizorkin spaces, all the spaces defined above do not depend on the choice of $\varphi$.
\item[ii)]
Previously we mentioned that our interest in these latter spaces was motivated by the work of Rivi\`ere and Struwe (see \cite{RS}) let us say a few words about this. In \cite{RS} the authors used the homogeneous Morrey space $L^{2,n-2}_1$ with norm
\begin{displaymath}
\vert \vert f \vert \vert_{L^{2,n-2}_1} ^2 = \sup_{x_0 \in \mathbb{R}^n} \sup_{r >0 } \Big ( \frac{1}{r^{n-2}}\int_{B-r(x_0)} \vert \nabla u \vert ^2 \Big).
\end{displaymath}
Note that $u \in L^{2,n-2}_1$ is equivalent to the fact that for all radii $r>0$ and all $x_0 \in \mathbb{R}^n$ we have the inequality
\begin{displaymath}
\vert \vert \nabla u \vert \vert_{L^2(B_r(x_0))} \leq C r^{(n-2)/p} = C r^{\frac{n}{2} -\frac{2}{2}}
\end{displaymath}
but this latter estimate is again equivalent to the fact that $\nabla u \in \mathcal{M}^n_2$.
Finally we remember that $\mathcal{M}^n_2 = \mathcal{N}^0_{n,2,2}$ (see for instance \cite{Maz2}) and note that $\nabla u \in \mathcal{N}^0_{n,2,2}$ is equivalent to $u \in \mathcal{N}^1_{n,2,2}$ since for all $s$ - even for the negative ones -  we have the equivalence $2^s \vert \vert u^s \vert \vert_{\mathcal{M}^n_2} \simeq \vert \vert (\nabla u)^s \vert \vert_{\mathcal{M}^n_2}$ because we always avoid the origin in the Fourier space and also near the origin work with annuli with radii $r \simeq 2^s$.
\end{itemize}
}
\end{Remark}

Before we continue, let us state a few facts concerning the spaces $B^s_{\mathcal{M}^p_q,r}$ which are interesting and important.

\begin{Lemma}\label{compl/scal}
\begin{itemize}
\item[i)]
The spaces $B^s_{\mathcal{M}^p_q,r}$ are complete for all possible choices of indices.
\item[ii)]
\begin{itemize}
\item[a)]
Let $s>0$, $1\leq q \leq p< \infty$, $1  \leq r \leq \infty$ and $\lambda > 0$. Then
\begin{displaymath}
\vert \vert f(\lambda \cdot ) \vert B^s_{\mathcal{M}^p_q,r}\vert \vert 
\leq C \lambda ^{-\frac{n}{p}}\sup \left \{ 1,\lambda \right \} ^s  \vert \vert f\vert B^s_{\mathcal{M}^p_q,r}\vert \vert.
\end{displaymath}
\item[b)]
Let $s=0$, $1\leq q \leq p< \infty$, $1 \leq r \leq \infty$ and $\lambda > 0$. Then
\begin{displaymath}
\vert \vert f(\lambda \cdot ) \vert B^s_{\mathcal{M}^p_q,r}\vert \vert 
\leq C \lambda ^{-\frac{n}{p}} (1 + \vert \log \lambda \vert )^{\alpha} \vert \vert f\vert B^s_{\mathcal{M}^p_q,r}\vert \vert
\end{displaymath}
\end{itemize}
where
\begin{displaymath}
\alpha=\frac{1}{r} \; \textrm{if} \; \lambda > 1 \; \; \textrm{and} \; \alpha=1-\frac{1}{r}=\frac{1}{r'} \; \textrm{if} \; 0< \lambda <1.
\end{displaymath}
\end{itemize}
\end{Lemma}
The first assertion is obtained by the same proof as the corresponding claim for the spaces $N^s_{p,q,r}$ in \cite{KY}.\\
The second fact is a variation of a well known proof given in \cite{B1}. \\

Furthermore we have the following embedding result which relates the spaces $B^0_{\mathcal{M}^p_q,r}$ to the Morrey spaces with the same indices respectively, similar for the spaces $N^0_{p,q,r}$.

\begin{Lemma}\label{inMorrey}
Let $1 < q \leq 2$, $1 < q \leq p < \infty$ and $r \leq q$. Then
\begin{displaymath}
B^0_{\mathcal{M}^p_q,r} \subset \mathcal{M}^p_q
\end{displaymath}
and
\begin{displaymath}
N^0_{p_q,r} \subset M^p_q.
\end{displaymath}
\end{Lemma}

From this result we immediately deduce the following corollary.

\begin{Corollary}
Let $1 < q \leq 2$, $1 < q \leq p < \infty$ and $r \leq q$ and assume that $f \in B^0_{\mathcal{M}^p_q,r}$ has compact support. Then $f \in L^q$.
\end{Corollary}

This holds because of the preceding lemma and the fact that for a bounded domain $\Omega$ we have the embedding $M^p_q(\Omega) \subset L^q(\Omega)$.

Similar to the result that $W^{1,p}=F^1_{p,2}$, $1<p<\infty$ we have the following lemma.

\begin{Lemma}\label{aequiv}
Assume that $f$ is a compactly supported distribution. Then, if $1 < q \leq 2$, $1 < q \leq p < \infty$ and $r \leq q$, the following two norms are equivalent
\begin{displaymath}
\vert \vert f \vert B^0_{\mathcal{M}^p_q,r} \vert \vert + \vert \vert \nabla f \vert B^0_{\mathcal{M}^p_q,r} \vert \vert
\end{displaymath}
\begin{displaymath}
\vert \vert f \vert B^1_{\mathcal{M}^p_q,r} \vert \vert.
\end{displaymath}
\end{Lemma} 
 
Moreover, also the fact that for a compactly supported distribution the homogeneous and the inhomogeneous Sobolev norms are equivalent, we have the following result.

\begin{Lemma} \label{aufKomp}
Let $1 < q \leq 2$, $1 < q \leq p < \infty$, $2 \leq p$, $r \leq q$ and $n\geq 3$.
Assume that the distribution $f$ has the following properties: $f$ has compact support and $\nabla f \in B^0_{\mathcal{M}^p_q,r}$. Then 
\begin{displaymath}
f \in B^1_{\mathcal{M}^p_q,r}.
\end{displaymath}
\end{Lemma}

As a by-product of our studies we have the following density result.

\begin{Lemma}\label{density}
Let $1\leq q \leq p < \infty$, $1\leq r \leq \infty$ and $s \in \mathbb{R}$. Then 
$O_M$ is dense in $N^s_{p,q,r}$ respectively in $\mathcal{N}^s_{p,q,r}$ and $B^s_{\mathcal{M}^p_q,r}$ where $O_M$ denotes the space of all $C^{\infty}$-functions such that 
$\forall \beta \in \mathbb{N}^n$ there exist constants $C_{\beta}>0$ and $m_{\beta } \in \mathbb{N}$ such that
\begin{displaymath}
\vert \partial^{\beta} f(x) \vert \leq C_{\beta} (1+\vert x \vert )^{m_{\beta}} \; \forall x \in \mathbb{R}^n.
\end{displaymath} 
Moreover, if $f \in N^s_{p,q,r}$ or $f\in B^s_{\mathcal{M}^q_pr,}$ with $s \geq 0$, $1\leq q \leq 2$ and $1 \leqq \leq p \leq \infty$ has compact support, it can be approximated by elements in $C^{\infty}_0$.
\end{Lemma}

Last, but not least we would like to mention a stability result which we will apply later on. 

\begin{Lemma}\label{stab}
Let $g \in B^0_{\mathcal{M}^n_2,2}$ and $f \in B^1_{\mathcal{M}^n_2,2,} \cap L^{\infty}$. Then
\begin{displaymath}
\vert \vert gf \vert B^0_{\mathcal{M}^n_2,2} \vert \vert \leq C \vert \vert g \vert B^0_{\mathcal{M}^n_2,2} \vert \vert (
\vert \vert f \vert B^1_{\mathcal{M}^n_2,2} \vert \vert + \vert \vert f \vert \vert_{\infty} ),
\end{displaymath}
i.e. $B^0_{\mathcal{M}^n_2,2}$ is stable under multiplication with a function in $B^1_{\mathcal{M}^n_2,2} \cap L^{\infty}$. 
\end{Lemma}

The proofs of lemma \ref{inMorrey}, \ref{aequiv}, \ref{aufKomp}, \ref{density} and \ref{stab} are given in the next section.\\

For further information about the Besov-Morrey spaces, see \cite{KY}, \cite{Maz1} and \cite{Maz2}.

\subsection{Spaces involving Choquet integrals}

In what follows, we will use a certain description of the pre-dual space of $\mathcal{M}^1$. Before we can state this assertion we have to introduce some function spaces involving the so-called Choquet integral. A general reference for this section is \cite{A} and the references given therein. \\
We start with the notion of Hausdorff capacity:

\begin{Definition}[Hausdorff capacity]
Let $E$ be a subset of $\mathbb{R}^n$ and let \\
$\left\{ B_j \right \}, \; j=1,2, \dots $ be a cover of $E$, i.e. $\left\{ B_j \right \}$ is a countable  collection of open balls $B_j$ with radius $r_j$ such that $E \subset \cup_j B_j$. Then we define the \textbf{Hausdorff capacity of $E$ of dimension $d$, $0 < d \leq n$} to be the following quantity
\begin{displaymath}
H^d_{\infty}(E)= \inf \sum_j r^d_j
\end{displaymath}
where the infimum is taken over all possible covers of $E$.
\end{Definition}

\begin{Remark}
\textnormal{
The name capacity may lead to confusion. Here we use this expression in the sense of N. Meyers. See \cite{Mey}, page 257.
}
\end{Remark}

Once we have this capacity, we can pass to the Choquet integral of $\phi \in C_0(\mathbb{R}^n)^+$:

\begin{Definition}[Choquet integral and $\mathbf{L^1(H^d_{\infty})}$]
Let $\phi \in C_0(\mathbb{R}^n)^+$. Then the \textbf{Choquet integral} of $\phi$ with respect to the Hausdorff capacity $H^d_{\infty}$ is defined to be the following Riemann integral:
\begin{displaymath}
\int \phi \; dH^d_{\infty} \equiv \int_0^{\infty} H^d_{\infty}\lbrack \phi > \lambda \rbrack \; d\lambda.
\end{displaymath}
The space $\mathbf{L^1(H^d_{\infty})}$ is now the completion of $C_0(\mathbb{R}^n)$ under the functional $\int \vert \phi \vert \; dH^d_{\infty}$.
\end{Definition}

Two important facts about $L^1(H^d_{\infty})$ are summarised below, again for instance see \cite{A} and also the references given there.

\begin{Remark}
\textnormal{
\begin{itemize}
\item
$L^1(H^d_{\infty})$ can also be characterised to be the space of all $H^d_{\infty}$-quasi continuous functions $\phi$ which satisfy $\int \vert \phi \vert \; dH^d_{\infty} < \infty$, i.e. for all $\varepsilon > 0$ there exists an open set $G$ such that $H^d_{\infty}[G] < \varepsilon$ and that $\phi$ restricted to the complement of $G$ is continuous there.
\item
One can show that $L^1(H^d_{\infty})$ is a quasi-Banach space with respect to the quasi-norm $\int \vert \phi \vert \; dH^d_{\infty}$. 
\end{itemize}
}
\end{Remark}

Now, we can state the duality result we mentioned earlier. A proof of this assertion is given in \cite{A}, but take care of the notation which differs from our notation!

\begin{Proposition}\label{dualityAdams}
We have $(L^1(H^d_{\infty}))^* =\mathcal{M}^{\frac{n}{n-d}}$ and in particular the estimate
\begin{displaymath}
\Big \vert \int u \; d\mu \Big \vert \leq \vert \vert u \vert \vert_{L^1(H^d_{\infty})} \vert \vert \mu \vert \vert_{\mathcal{M}^{\frac{n}{n-d}}}
\end{displaymath}
holds and
\begin{displaymath}
\vert \vert \mu \vert \vert_{(L^1(H^d_{\infty}))^*}=\sup_{\vert\vert u \vert \vert _{L^1(H^d_{\infty})} \leq 1} \Big \vert \int u \; d\mu \Big \vert 
\simeq \vert \vert \mu \vert \vert _{\mathcal{M}^{\frac{n}{n-d}}}.
\end{displaymath}
\end{Proposition}

Note that in order to show that a certain function belongs to $\mathcal{M}^{\frac{n}{n-d}}$, it is enough to show that it defines a linear functional on $L^1(H^d_{\infty})$, i. e. that $\sup_{\vert\vert u \vert \vert _{L^1(H^d_{\infty})} \leq 1} \Big \vert \int u \; d\mu \Big \vert < \infty$. This does not require that $L^1(H^d_{\infty})$ is a Banach space and is quite different from the case when you use the dual characterisation of a norm in order to show that a certain distribution belongs to a certain space. 

\begin{Remark}
\textnormal{
The above proposition is just a special case of a more general result which involves also spaces $L^p(H^d_{\infty})$, see for instance \cite{A2}.
}
\end{Remark}

Before ending this section we will state some useful remarks for later applications.

\begin{Remark}\label{ChoquetinTempdist}
\textnormal{
\begin{itemize}
\item
Observe that $\mathcal{M}^p \subset \mathcal{S}'$ (in particular for $p=\frac{n}{n-d}$). In order to verify this, note that $\mathcal{M}^p \subset N^0_{p,1,\infty} \subset \mathcal{S}'$: Let $\mu \in \mathcal{M}^p$ and let as usual $\varphi \in \Phi(\mathbb{R}^n)$ then we have
\begin{eqnarray}
\vert \vert \mu \vert N^0_{p,1,\infty} \vert \vert &=& \sup_{k \in \mathbb{N}} \vert \vert \check \varphi_k \ast \mu \vert M^p_1 \vert \vert \nonumber \\
&=&\sup_{k \in \mathbb{N}} \vert \vert \check \varphi_k \ast \mu \vert M^p \vert \vert \nonumber \\
&&\textrm{note that $\check \varphi_k \ast \mu \in C^{\infty} \subset L^1_{loc}$ since $\mu \in \mathcal{D}'$} \nonumber \\
&&\textrm{and $\check \varphi_k \ast \mu$ can be seen as a measure} \nonumber \\
&\leq& \sup_{k \in \mathcal{N}} \vert \vert \check \varphi_k \vert \vert_1 \vert \vert \mu \vert M^p \vert \vert \nonumber \\
&&\textrm{because of \cite{KY}, lemma 1.8} \nonumber \\
&\leq& C \vert \vert \mu \vert \mathcal{M}^p \vert \vert \nonumber \\
&<& \infty \nonumber \\
&&\textrm{according to our hypothesis.} \nonumber
\end{eqnarray} 
Once we have this, we apply the continuous embedding of $N^0_{p,1,\infty}$ into $\mathcal{S}'$ (see e.g. \cite{Maz2}) and conclude that actually $\mathcal{M}^p \subset \mathcal{S}'$. \\
Note also that $\mathcal{S} \subset L^1(H^d_{\infty})$
\item
Using the duality asserted above, we can show that $L^1(H^d_{\infty}) \subset \mathcal{S}'$:
We start with $f \in C^{\infty}_0(\mathbb{R}^n)$. Since $f \in L^{\infty}$ it is easy to check that $f \in M^p_q \; , 1 \leq q \leq p < \infty,$ with $\vert \vert f \vert M^p_q \vert \vert = \vert \vert f \vert \vert_{\infty}$. Moreover, $f$ even belongs to $\mathcal{M}^p_q$. In order to establish this, it remains to show that there is a constant $C$, independent on $f$, such that $\forall x \in \mathbb{R}^n$ and for $1\leq r$
\begin{displaymath}
\vert \vert f \vert \vert_{L^1(B_r(x))} \leq C r^{\frac{n}{q}-\frac{n}{p}}.
\end{displaymath}
In fact, it holds $\forall x \in \mathbb{R}^n$ and $\forall r \geq 1$
\begin{eqnarray}
\vert \vert f \vert \vert_{L^1(B_r(x))} & \leq& \vert \vert f \vert \vert_1 \nonumber \\
&\leq& \vert \vert f \vert \vert_1 r^{\frac{n}{q}-\frac{n}{p}}\nonumber \\
&&\textrm{since due to the choice of $p$ and $q$ we have} \nonumber \\
&&\frac{n}{q} - \frac{n}{p}\geq 0. \nonumber
\end{eqnarray}
If we put together all these information we find 
\begin{displaymath}
\vert \vert f \vert \mathcal{M}^p_q \vert \vert \leq \vert \vert f \vert \vert_{\infty}+ \vert \vert f \vert \vert_1.
\end{displaymath}
Now, recall that the duality between $L^1(H^d_{\infty})$ and $\mathcal{M}^{\frac{n}{n-d}}$ is given by
\begin{displaymath}
<\mu,u>_{(L^1(H^d_{\infty}))^*=\mathcal{M}^{\frac{n}{n-d}}, L^1(H^d_{\infty})}=\int u \; d \mu  
\end{displaymath}
where $u \in L^1(H^d_{\infty})$ and $\mu \in \mathcal{M}^{\frac{n}{n-d}}$.\\
In a next step we define the action of $u \in L^1(H _{\infty})$ on $f \in C^{\infty}_0$ as follows
\begin{displaymath}
< u, f>_{\mathcal{D}',C^{\infty}_0}:=< f,u>_{\mathcal{M}^{\frac{n}{n-d}}, L^1(H^d_{\infty})}.
\end{displaymath}
Last, but not least, we observe that for $\varphi \in \mathcal{S}$ we have
\begin{displaymath}
\vert \vert \varphi \vert \vert_{\infty} + \vert \vert \varphi \vert \vert_1 \leq C(n) \vert \vert \varphi \vert \vert _{\mathcal{S}}.
\end{displaymath}
This finally leads to the conclusion that in fact, $L^1(H^d_{\infty}) \subset \mathcal{S}'$.
\end{itemize}
}
\end{Remark}

This last remark enables us to use the above introduced $L^1(H^d_{\infty})$-quasi norm to construct - in analogy to the case of Besov- or Besov-Morrey-spaces - a new space of functions.

\begin{Definition}[Besov-Choquet spaces]
Let $\varphi \in \Phi (\mathbb{R}^n)$.\\
We say that $f\in \mathcal{S}'$ belongs to $\mathbf{B^0_{L^1(H^d_{\infty}),\infty}}$ if
$\exists \left\{f_k(x)\right\} _{k=0}^{\infty} \subset L^1(H^d_{\infty})$ such that the following holds
\begin{displaymath}
f=\sum_{k=0}^{\infty} \mathcal{F}^{-1}\varphi _k\mathcal{F}f_k  \; \textrm{in} \; \mathcal{S}'(\mathbb{R}^n)
\end{displaymath}
and 
\begin{displaymath}
\sup_k \vert \vert f_k  \vert  L^1(H^d_{\infty}) \vert \vert < \infty.
\end{displaymath}
Moreover we set
\begin{displaymath}
\vert \vert f \vert B^0_{L^1(H^d_{\infty}),\infty} \vert \vert = 
\inf \sup_k \vert \vert f_k  \vert L^1(H^d_{\infty})\vert\vert
\end{displaymath} 
where the infimum is taken over all admissible representations of $f$. \\
Moreover, we denote by $\mathbf{b^0_{L^1(H^d_{\infty}),\infty}}$ the closure of $\mathcal{S}$ under the construction explained above.
\end {Definition}

\begin{Remark}
\textnormal{
In complete analogy to the construction of the Besov spaces (respectively the Besov-Morrey-spaces) one could also construct new spaces if we replace the Lebesgue $L^p$-norms (respectively the Morrey-norms) by $L^p(H^d_{\infty})$-quasi-norms.
}
\end{Remark}

\section{Proofs}

\subsection{Some preliminary remarks}

In what follows we set
\begin{displaymath}
\mathbf{f^j(x)}=\mathcal{F}^{-1}( \varphi _j\mathcal{F}f) (x)
\end{displaymath}
where $\varphi =\left\{\varphi _j(x)\right\}_{j=0}^{\infty} \in \Phi (\mathbb{R}^n)$.

Recall that once we can control the \textbf{paraproducts} $\pi_1(f,g)=\sum_{k=2}^{\infty}\sum_{l=0}^{k-2}f^lg^k$, $\pi_2(f,g)=\sum_{k=0}^{\infty}\sum_{l=k-1}^{k+1}f^lg^k $ and $\pi_3(f,g)=\sum_{l=2}^{\infty}\sum_{k=0}^{l-2}f^lg^k$ ($f^i=0$ if $i\leq -1$ and similarly for $g$) we are also able to control the product $fg$ (see e.g. \cite {RS}). \\
Since in the sequel we want to take into account cancellation phenomena, we will analysis
\begin{equation}\label{gleichung1}
\pi _1(a_x,b_y) \; , \; \pi_1(a_y,b_x) \; , \pi _3(a_x,b_y) \; , \pi_3(a_y,b_x) \; \textrm{and} 
\; \sum_{s=0}^{\infty}\sum_{t=s-1}^{s+1}a_x^tb_y^s -a_y^tb_x^s. 
\end{equation}
Last but not least, remember that
\begin{displaymath}
supp \; \mathcal{F}\Big( \sum_{i=0}^{l-2}a_x^i b_y^l \Big) \subset 
\left\{ \xi \; \vert \; 2^{l-3}\leq \vert \xi \vert \leq 2^{l+3} \right\}\textrm{for} \; l\geq 2.
\end{displaymath}
and
\begin{displaymath}
supp \; \mathcal{F}\Big( \sum_{i=l-1}^{l+1}a_x^i b_y^l \Big) \subset 
\left\{ \xi \; \vert \; \vert \xi \vert \leq 5\cdot 2^{l} \right\}\textrm{for} \; l\geq 0.
\end{displaymath} 

\subsection{Proof of theorem \ref{regndim} i)}

The proof of this assertion is split into several parts: In a first step we show that $\pi_1(a_x,b_y), \; \pi_3(a_x,b_y), \; \pi_3(a_y,b_x)$ and $\pi_1(a_y,b_x) \in B^{-1}_{\infty,1}$ and $\sum_{s=0}^{\infty}\sum_{t=s-1}^{s+1}a_x^tb_y^s -a_y^tb_x^s \in B^{-2}_{\infty,1}$. Once we have this we show in a second step that under this hypothesis the solution $u$ of
\begin{equation}
-\Delta u =f \; \textrm{where} f \in B^{-2}_{\infty,1} \nonumber
\end{equation}
is continuous.\\

\textbf{\textit{Claim: $\pi_1(a_x,b_y) \in B^{-2}_{\infty,1}$}}\\

Our hypotheses together with \cite{KY}, theorem 2.5, ensures us that $a_x, \; b_y \in B^{-1}_{\infty,2}$. Next, due to \cite{RS}, proposition 1, chapter 2.3.2,  it is enough to prove that
\begin{equation}
\vert \vert2^{-2j}c_j \vert l^1(L^{\infty}) \vert \vert < \infty \nonumber
\end{equation}
where as before $c_j := \sum_{t=0}^{k-2}a_x^tb_y^j$. \\
We actually have
\begin{eqnarray}
\vert \vert2^{-2j}c_j \vert l^1(L^{\infty}) \vert \vert &=& \sum_{j=0}^{\infty} 2^{-2j} \vert \vert \sum_{t=0}^{^j-2}a_x^tb_y^j \vert \vert_{\infty} \nonumber \\
&\leq& \sum_{j=0}^{\infty} 2^{-2j} \vert \vert \sum_{t=0}^{j-2}a_x^t \vert \vert_{\infty} \vert \vert b_y^j \vert \vert_{\infty} \nonumber \\
&=& \sum_{j=0}^{\infty} 2^{-j} \vert \vert \sum_{t=0}^{j-2}a_x^t \vert \vert_{\infty} 2^{-j}\vert \vert b_y^j \vert \vert_{\infty} \nonumber \\
&\leq& \Big( \sum_{j=0}^{\infty} 2^{-2j} \vert \vert \sum_{t=0}^{j-2}a_x^t \vert \vert_{\infty}^2 \Big)^{\frac{1}{2}} \Big( \sum_{j=0}^{\infty} 2^{-2j}\vert \vert b_y^j \vert \vert_{\infty}^2 \Big)^{\frac{1}{2}} \nonumber \\
&&\textrm{due to H\"older's inequality} \nonumber \\
&=& \vert \vert 2^{-j} \vert \sum_{t=0}^{j-2} a_x^t \vert \; \vert l^2(L^{\infty}) \vert \vert \vert \vert b_y \vert B^{-1}_{\infty} \vert \vert \nonumber \\
&\leq& C \vert \vert 2^{-j} \vert \sum_{t=0}^{j} a_x^t \vert \; \vert l^2(L^{\infty}) \vert \vert \vert \vert b_y \vert B^{-1}_{\infty} \vert \vert \nonumber \\
&\leq& C \vert \vert a_x \vert B^{-1}_ {\infty,2} \vert \vert \vert \vert b_y \vert B^{-1}_{\infty} \vert \vert \; \textrm{because of \cite{RS}, first lemma in chapter 4.4.2 }\nonumber \\
&<&\infty \; \textrm{thanks to our hypothesis.} \nonumber
\end{eqnarray}

This shows that in fact $\pi_1(a_x,b_y) \in B^{-2}_{\infty,1}$. Similarly one proves that also $\pi_1(a_y,b_x), \; \pi_3(a_x,b_y) \; \textrm{and} \; \pi_1(a_y,b_x)$ belong to the same space.

In remains to analyse the contribution where the frequencies are comparable. This is our next goal. \\

\textbf{\textit{Analysis of $\sum_{s=0}^{\infty}\sum_{t=s-1}^{s+1}a_x^tb_y^s -a_y^tb_x^s$}}\\

In stead of first applying the embedding result of Kozono/Yamazaki which embeds Morrey-Besov spaces into Besov spaces and then analysing a certain quantity, we invert the order of these steps in order to estimate $\sum_{s=0}^{\infty}\sum_{t=s-1}^{s+1}a_x^tb_y^s -a_y^tb_x^s$.\\
We will use the following result concerning predual spaces of Morrey spaces .

\begin{Proposition}\label{Duality}
The dual space of $b^0_{L^1(H^{n-2}_\infty),\infty}$ is the space $B^0_{\mathcal{M}^{\frac{n}{2}}_1,1}$. 
\end{Proposition}

\begin{Remark}
\textnormal{
The above result has the same flavour as (see for instance \cite{RS})
\begin{displaymath}
(b^0_{\infty,\infty})^*=B^0_{1,1}
\end{displaymath}
}
\end{Remark}

\textit{Proof of proposition \ref{Duality}:}\\
We have to show the two inclusion relations. \\
We start with $(b^0_{L^1(H^{n-2}_\infty),\infty})^* \supset B^0_{\mathcal{M}^{\frac{n}{2}}_1,1}$: \\
Assume that $f \in B^0_{\mathcal{M}^{\frac{n}{2}}_1,1}\subset N^0_{\frac{n}{2},1,1} \subset \mathcal{S}'$ and assume that $\psi \in b^0_{L^1(H^{n-2}_\infty),\infty}$. By density we may assume that $\psi \in \mathcal{S}$. We have to show that $f \in (b^0_{L^1(H^{n-2}_\infty),\infty})^*$. To this end let $\sum_{k=0}^{\infty} \check \varphi_k \ast \psi_k$ be a representation of $\psi$ with 
\begin{displaymath}
\sup_{k} \vert \vert \psi_k \vert \vert _{L^1(H^{n-2}_{\infty})} \leq 2 \vert \vert \psi \vert b^0_{L^1(H^{n-2}_\infty),\infty} \vert \vert.
\end{displaymath}
Note that in our case - as a tempered distribution - $f$ acts on $\psi$ and we estimate
\begin{eqnarray}
\vert f(\psi)\vert &=& \vert f(\sum_{k\geq 0} \check \varphi_k \ast\psi_k)\vert \nonumber \\
&=& \Big\vert f\Big( \sum_{k=o}^{\infty} \mathcal{F}^{-1}(\varphi_k\mathcal{F}\psi_k)\Big)\Big\vert \nonumber \\
&=&\Big \vert \sum_{k=0}^{\infty}f\Big(\mathcal{F}^{-1}(\varphi_k \mathcal{F}\psi_k)\Big) \Big \vert = 
\Big \vert \sum_{k=0}^{\infty}\int f\mathcal{F}^{-1}(\varphi_k \mathcal{F}\psi_k\Big) \Big \vert\nonumber \\
&=& \Big \vert \sum_{k=0}^{\infty}\psi_k\mathcal{F}(\varphi_k\mathcal{F}^{-1}f) \Big \vert =
\Big \vert \sum_{k=0}^{\infty}\int\psi_k \; d f  \Big \vert \nonumber \\
&&\textrm{where $df= \mathcal{F}(\varphi_k \mathcal{F}^{-1}f)\;d\lambda$ with $\lambda$ the Lebesgue measure} \nonumber \\
&\leq& \sum_{k=0}^{\infty} \vert \psi_k \mathcal{F}(\varphi_k\mathcal{F}^{-1}f) \vert \nonumber \nonumber \\
&\leq& \sup_{k\geq 0}\vert \vert\psi_k \vert \vert_{L^1(H^{n-2}_{\infty})} \sum_{k=o}^{\infty} \vert \vert \mathcal{F}(\varphi_k\mathcal{F}^{-1}f) \vert \vert _{\mathcal{M}^{\frac{n}{2}}} \; \textrm{recall proposition \ref{dualityAdams}}\nonumber \\
&=& \sup_{k\geq 0}\vert \vert\psi_k \vert \vert_{L^1(H^{n-2}_{\infty})} \sum_{k=o}^{\infty} \vert \vert \mathcal{F}(\varphi_k\mathcal{F}^{-1}f )\vert \vert _{\mathcal{M}^{\frac{n}{2}}_1}\nonumber \\
&&\textrm{cf. also remark \ref{ChoquetinTempdist}} \nonumber \\
&\leq& C \sup_{k\geq 0}\vert \vert\psi_k \vert \vert_{L^1(H^{n-2}_{\infty})} \sum_{k=o}^{\infty} \vert \vert \mathcal{F}^{-1}(\varphi_k\mathcal{F}f )\vert \vert _{\mathcal{M}^{\frac{n}{2}}_1}\nonumber \\
&\leq& C \vert \vert \psi \vert b^0_{L^1(H^{n-2}_{\infty}), \infty} \vert \vert \; \;\vert \vert f \vert B^0_{\mathcal{M}^{\frac{n}{2}}_1,1} \vert \vert \nonumber \\
&<& \infty \nonumber \\
&&\textrm{thanks to our assumptions.} \nonumber 
\end{eqnarray}
Now we show the other inclusion, $(b^0_{L^1(H^{n-2}_\infty),\infty})^* \subset B^0_{\mathcal{M}^{\frac{n}{2}}_1,1}$: \\
We start with $f \in (b^0_{L^1(H^{n-2}_\infty),\infty})^*$ and we have to show that $f$ belongs also to $B^0_{\mathcal{M}^{\frac{n}{2}}_1,1}$:
First of all, note that $f$ gives also rise to elements of $(L^1(H^{n-2}_{\infty}))^*$ as follows: Each $\psi \in b^0_{L^1(H^{n-2}_{\infty}),\infty}$ can be seen as a sequence $\left\{\psi_k\right\}_{k=0}^{\infty} \subset L^1(H^{n-2}_{\infty})$, and of course $\check \varphi_k \ast \psi_k \in b^0_{L^1(H^{n-2}_{\infty}),\infty} \; \forall k \in \mathbb{N}$. Moreover, for each $k \in \mathbb{N}$ we have - again by density of $\mathcal{S}$ -
\begin{eqnarray}
f(\delta_{kj} (\check \varphi_j\ast \psi_j))&=&<f,\delta_{kj}\psi>_{(b^0_{L^1(H^{n-2}_\infty),\infty})^*,b^0_{L^1(H^{n-2}_\infty),\infty}}\nonumber \\
&=& <f,\check \varphi_k \ast \psi_k >_{(b^0_{L^1(H^{n-2}_\infty),\infty})^*,b^0_{L^1(H^{n-2}_\infty),\infty}} \nonumber \\
&=&<f,\check \varphi_k \ast \psi_k >_{\mathcal{S}',\mathcal{S}}\nonumber \\
&=&<f,\mathcal{F}^{-1}( \varphi_k \mathcal{F} \psi_k )>_{\mathcal{S}',\mathcal{S}}\nonumber \\
&=&<\mathcal{F}(\varphi_k \mathcal{F}^{-1}f),\psi_k >_{\mathcal{S}',\mathcal{S}}\nonumber \\
&=&<\mathcal{F}(\varphi_k \mathcal{F}^{-1}f),\psi_k >_{\mathcal{M}^{\frac{n}{2}},L^1(H^{n-2}_{\infty})}.\nonumber 
\end{eqnarray}
Next we will construct a special element of $b^0_{L^1(H^{n-2}_\infty),\infty}$: \\
Let $0<\varepsilon$ small.\\
We choose $\psi_k$ such that 
\begin{itemize}
\item
$\psi_k  \in \mathcal{S}$: Remember that we have density!
\item
$\vert \vert \psi_k \vert \vert_{L^1(H^{n-2}_{\infty})} \leq 1$
\item
$0< \; \; <\mathcal{F}(\varphi_k \mathcal{F}^{-1}f),\psi_k >_{\mathcal{M}^{\frac{n}{2}},L^1(H^{n-2}_{\infty})}$
\item
\begin{eqnarray}
<\mathcal{F}(\varphi_k \mathcal{F}^{-1}f),\psi_k >_{\mathcal{M}^{\frac{n}{2}},L^1(H^{n-2}_{\infty})} &\geq& \vert \vert \mathcal{F}(\varphi_k \mathcal{F}^{-1}f)\vert \vert_{\mathcal{M}^{\frac{n}{2}}}- \varepsilon 2^{-k}\nonumber \\
&=&  \vert \vert \mathcal{F}(\varphi_k \mathcal{F}^{-1}f)\vert \vert_{(L^1(H^{n-2}_{\infty}))^*} -\varepsilon 2^{-k}\nonumber\\
&=& \sup_{\substack{u \in L^1(H^{n-2}_{\infty})\\ \vert \vert u \vert \vert _{L^1(H^{n-2}_{\infty})}\leq 1}} \vert <\mathcal{F}(\varphi_k \mathcal{F}^{-1}f),u> \vert -\varepsilon 2^{-k}. \nonumber
\end{eqnarray}
\end{itemize}
Note that like that $\psi=\sum_{k=0}^{\infty} \check \varphi_k \ast \psi_k \in b^0_{L^1(H^{n-2}_\infty),\infty}$ with 
\begin{displaymath}
\vert \vert \psi \vert b^0_{L^1(H^{n-2}_\infty),\infty} \vert \vert \leq 1.
\end{displaymath}
If we put now all this together we find - recall that $f$ acts linearly! - 
\begin{eqnarray}
\sum_{k=0}^{\infty} \vert \vert f^k \vert \vert_{\mathcal{M}^{\frac{n}{2}}_1} &=& 
\sum_{k=0}^{\infty} \vert \vert \mathcal{F}^{-1}(\varphi_k \mathcal{F}f) \vert \vert_{\mathcal{M}^{\frac{n}{2}}_1}\nonumber \\
&=& C \sum_{k=0}^{\infty} \vert \vert \mathcal{F}(\varphi_k \mathcal{F}^{-1}f) \vert \vert_{\mathcal{M}^{\frac{n}{2}}_1} \nonumber \\
&\leq& 2 \varepsilon + f(\psi) \nonumber \\
&&\textrm{$\psi$ as constructed above}\nonumber \\
&\leq& 2 \varepsilon  + \vert \vert f \vert  (b^0_{L^1(H^{n-2}_\infty),\infty})^* \vert \vert \; \; \vert \vert \psi \vert b^0_{L^1(H^{n-2}_\infty),\infty}\vert \vert \nonumber \\
&\leq& 2\varepsilon + \vert \vert f \vert  (b^0_{L^1(H^{n-2}_\infty),\infty})^* \vert \vert. \nonumber
\end{eqnarray}
Since this holds for all $0<\varepsilon$ we let $\varepsilon$ tend to zero and get the desired inclusion.\\
All together we established the duality result we claimed above. 
\begin{flushright}
$\Box$
\end{flushright}

What concerns the next lemma, recall that $\mathcal{S}$ is dense in $b^0_{L^1(H^{n-2}_\infty),\infty}$:

\begin{Lemma}\label{ndimmainlemma}
Let $\phi \in \Phi(\mathbb{R}^n)$ and assume that $\psi  \in \mathcal{S} \cap L^1(H^{n-2}_{\infty})$ with representation $\left\{\psi_k \right\}_{k=0}^{\infty}$, i.e.
$\sum_{k=0}^{\infty} \check \varphi_k \ast \psi_k=\psi$, such that
\begin{displaymath}
\sup_{k} \vert \vert \psi_k \vert \vert _{L^1(H^{n-2}_{\infty})} \leq 2 \vert \vert \psi \vert b^0_{L^1(H^{n-2}_\infty),\infty} \vert \vert.
\end{displaymath}
Then
\begin{eqnarray}
\Big\vert\Big \vert \frac{\partial}{\partial x}\check \varphi_k \ast \psi_k  \Big\vert \Big\vert_{L^1(H^{n-2}_{\infty})}
&=& \Big\vert\Big \vert \frac{\partial}{\partial x} (\check{\varphi_k} \ast \psi_k) \Big\vert \Big\vert_{L^1(H^{n-2}_{\infty})}\nonumber \\
&\leq& C2^s \vert \vert \psi_k \vert \vert_{L^1(H^{n-2}_{\infty})}
\leq C2^s \vert \vert \psi \vert b^0_{L^1(H^{n-2}_{\infty}),\infty}\vert \vert. \nonumber
\end{eqnarray}
\end{Lemma}

\textit{Proof:}\\
For the proof of this lemma, we need the fact that if $f(x)\geq 0$ is lower semi-continuous on $\mathbb{R}^n$ then
\begin{displaymath}
\vert \vert f \vert \vert_{L^1(H^{d}_{\infty})} = \int f \; dH^d_{\infty} \sim \sup \left\{ \int f \; d\mu \vert \mu \in \mathcal{M}^{\frac{n}{n-d}}_+ \; \textrm{and} \; \vert \vert \mu \vert \vert_{\mathcal{M}^{\frac{n}{n-d}}} \leq 1 
\right\}.
\end{displaymath}
see Adams \cite{A}.\\

It holds
\begin{eqnarray}
\Big\vert \Big\vert \frac{\partial}{\partial x}\check \varphi_k \ast \psi_k  \Big\vert \Big \vert_{L^1(H^{n-2}_{\infty})}
&\leq& \vert \vert \; \vert \frac{\partial}{\partial x} \check\varphi_k \vert \ast \vert\psi_k \vert \; \vert \vert_{L^1(H^{n-2}_{\infty})}\nonumber \\
&\leq& C \sup _{\substack{\mu  \in \mathcal{M}^{\frac{n}{2}}_+\\ \vert \vert \mu \vert \vert _{\mathcal{M}^{\frac{n}{2}}}\leq 1}}\left\{ \int  \vert \frac{\partial}{\partial x} \check\varphi_k \vert \ast \vert\psi_k \vert\; d\mu \right\}\nonumber \\
&=& C \sup _{\substack{\mu  \in \mathcal{M}^{\frac{n}{2}}_+\\ \vert \vert \mu \vert \vert _{\mathcal{M}^{\frac{n}{2}}}\leq 1}}\left\{ \int  \int \vert \frac{\partial}{\partial x} \check\varphi_k \vert(x-y)  \vert\psi_k \vert(y)\; d\lambda(y)d\mu(x) \right\}\nonumber \\
&\leq& C \sup _{\substack{\mu  \in \mathcal{M}^{\frac{n}{2}}_+\\ \vert \vert \mu \vert \vert _{\mathcal{M}^{\frac{n}{2}}}\leq 1}}\left\{ \int  \vert \psi_k \vert (y) \int \vert \frac{\partial}{\partial x} \check\varphi_k \vert(x-y)  \;  d\mu(x)d\lambda(y) \right\}\nonumber \\
&&\textrm{by Tonelli's theorem} \nonumber \\
&=&C\sup _{\substack{\mu  \in \mathcal{M}^{\frac{n}{2}}_+\\ \vert \vert \mu \vert \vert _{\mathcal{M}^{\frac{n}{2}}}\leq 1}}\left\{ \int  \vert \psi_k \vert (y) \int \vert \frac{\partial}{\partial x} \check\varphi_k \vert(y-x)  \;  d\mu(x)d\lambda(y) \right\}\nonumber \\
&&\textrm{note that $\varphi_k$ can be chose radial} \nonumber \\
&&\textrm{which implies that $\check \varphi_k$ and $\frac{\partial}{\partial x} \check\varphi_k$ are radial} \nonumber \\ 
&&\textrm{see e.g. \cite{SW}} \nonumber \\
&=& C\sup _{\substack{\mu  \in \mathcal{M}^{\frac{n}{2}}_+\\ \vert \vert \mu \vert \vert _{\mathcal{M}^{\frac{n}{2}}}\leq 1}}\left\{ \int  \vert \psi_k \vert (y) \frac{\partial}{\partial x} \check\varphi_k \vert(y-x) \ast \mu(y) \; d\lambda(y) \right\}\nonumber 
\end{eqnarray}
and we continue
\begin{eqnarray}
\Big\vert \Big\vert \frac{\partial}{\partial x}\check \varphi_k \ast \psi_k  \Big\vert \Big \vert_{L^1(H^{n-2}_{\infty})}
&\leq & C\sup _{\substack{\mu  \in \mathcal{M}^{\frac{n}{2}}_+\\ \vert \vert \mu \vert \vert _{\mathcal{M}^{\frac{n}{2}}}\leq 1}}\left\{ \int  \vert \psi_k \vert (y)\; d\nu(y) \right\}\nonumber \\
&&\textrm{where $\nu:=\frac{\partial}{\partial x} \check\varphi_k \lambda \ast \mu$}\nonumber \\
&\leq& C\sup _{\substack{\mu  \in \mathcal{M}^{\frac{n}{2}}_+\\ \vert \vert \mu \vert \vert _{\mathcal{M}^{\frac{n}{2}}}\leq 1}}\left\{ \vert \vert \psi_k \vert \vert_{L^1(H^{n-2}_{\infty})} \; \vert \vert \frac{\partial}{\partial x} \check\varphi_k \lambda \ast \mu \vert \vert_{\mathcal{M}^{\frac{n}{2}}}\right\}\nonumber \\
&\leq& C\sup _{\substack{\mu  \in \mathcal{M}^{\frac{n}{2}}_+\\ \vert \vert \mu \vert \vert _{\mathcal{M}^{\frac{n}{2}}}\leq 1}}\left\{ \vert \vert \psi_k \vert \vert_{L^1(H^{n-2}_{\infty})} \; \vert \vert \frac{\partial}{\partial x} \check\varphi_k \vert \vert_{L^1} \; \vert \vert \mu \vert \vert_{\mathcal{M}^{\frac{n}{2}}}\right\}\nonumber \\
&&\textrm{by \cite{KY}, lemma 1.8} \nonumber \\
&\leq& C\vert \vert \psi_k \vert \vert_{L^1(H^{n-2}_{\infty})} \; \vert \vert \frac{\partial}{\partial x} \check\varphi_k \vert \vert_{L^1} \nonumber \\
&\leq& C2^k \vert \vert \psi_k \vert \vert_{L^1(H^{n-2}_{\infty})}\nonumber \\
&\leq& C 2^k \vert \vert \psi \vert b^0_{L^1(H^{n-2}_\infty),\infty} \vert \vert  \nonumber
\end{eqnarray}
what we had to prove.
\begin{flushright}
$\Box$
\end{flushright}

The next lemma is a technical one:

\begin{Lemma} \label{PartSum}
Let $a \; \textrm{and} \; b \; \textrm{belong to} \; C^{\infty}_0(\mathbb{R}^n)$, $t=s+j$ where $j \in \left\{-1,0,1 \right \}$ and 
$\psi$ with representation $\left\{\psi_k\right\}_{k=0}^{\infty}$, i.e. $\sum_{k=0}^{\infty}\check \varphi_k \ast \psi_k = \psi$, such that
\begin{displaymath}
\sup_{k} \vert \vert \psi_k \vert \vert _{L^1(H^{n-2}_{\infty})} \leq 2 \vert \vert \psi \vert b^0_{L^1(H^{n-2}_\infty),\infty} \vert \vert \leq 2
\end{displaymath}
Then
\begin{eqnarray}
&&\int_{\mathbb{R}^n} \frac{\partial}{\partial x} (a^t b_y^s)\psi - \frac{\partial}{\partial y}(a^tb_x^s)\psi \nonumber \\
&&=\int_{\mathbb{R}^n} \frac{\partial}{\partial x} \Big(a^t b_y^s\Big)\Big(\sum_{k=0}^{s+3}\mathcal{F}^{-1}(\varphi_k\mathcal{F}\psi_k)\Big)
-\frac{\partial}{\partial y}\Big(a^tb_x^s\Big)\Big(\sum_{k=0}^{s+3}\mathcal{F}^{-1}(\varphi_k\mathcal{F}\psi_k)\Big). \nonumber
\end{eqnarray}
\end{Lemma}

\textit{Proof:}\\
First of all, note that $h \in \mathcal{S}'$ and $a^tb_y^s$ and $a^tb_x^s$ belong to $\mathcal{S}$ independently of the choices of $s$ and $t$. \\
We now calculate
\begin{eqnarray}
\int_{\mathbb{R}^n}\frac{\partial}{\partial x} (a^tb_y^s )\psi -\frac{\partial}{\partial y}(a^tb_x^s)\psi&=&
\int_{\mathbb{R}^n}\frac{\partial}{\partial x} (a^tb_y^s )\psi -\int_{\mathbb{R}^n}\frac{\partial}{\partial y}(a^tb_x^s)\psi \nonumber \\
&=& \int_{\mathbb{R}^n}\frac{\partial}{\partial x} (a^tb_y^s )\sum_{k=0}^{\infty} \mathcal{F}^{-1}(\varphi_k \mathcal{F}\psi_k)\nonumber \\
&& -\int_{\mathbb{R}^n}\frac{\partial}{\partial y}(a^tb_x^s)\sum_{k=0}^{\infty} \mathcal{F}^{-1}(\varphi_k \mathcal{F}\psi_k) \nonumber 
\end{eqnarray}
so
\begin{eqnarray}
\int_{\mathbb{R}^n}\frac{\partial}{\partial x} (a^tb_y^s )\psi -\frac{\partial}{\partial y}(a^tb_x^s)\psi
&=& \int_{\mathbb{R}^n}\frac{\partial}{\partial x} (a^tb_y^s )
\Big\lbrack \sum_{k=0}^{s+3} \mathcal{F}^{-1}(\varphi_k \mathcal{F}\psi_k) +\sum_{k=s+4}^{\infty}\mathcal{F}^{-1}(\varphi_k \mathcal{F}\psi_k)\Big\rbrack\nonumber \\
&& -\int_{\mathbb{R}^n}\frac{\partial}{\partial y}(a^tb_x^s)
\Big\lbrack\sum_{k=0}^{s+4} \mathcal{F}^{-1}(\varphi_k \mathcal{F}\psi_k) +\sum_{k=s+4}^{\infty}\mathcal{F}^{-1}(\varphi_k \mathcal{F}\psi_k)\Big\rbrack.\nonumber 
\end{eqnarray}
These calculations show that we have to prove that
\begin{eqnarray}
&&\int_{\mathbb{R}^n}\frac{\partial}{\partial x} (a^tb_y^s )
\sum_{k=s+4}^{\infty}\mathcal{F}^{-1}(\varphi_k \mathcal{F}\psi_k)=0 \nonumber \\
&&\textrm{and} \nonumber \\
&& \int_{\mathbb{R}^n}\frac{\partial}{\partial y} (a^tb_x^s )
\sum_{k=s+4}^{\infty}\mathcal{F}^{-1}(\varphi_k \mathcal{F}\psi_k)=0. \nonumber 
\end{eqnarray}
In what follows, we will only discuss the first integral because the second one can be analysed in exactly the same way. \\
So from now on we look at
\begin{displaymath}
\int_{\mathbb{R}^n}\frac{\partial}{\partial x} (a^tb_y^s )
\sum_{k=s+4}^{\infty}\mathcal{F}^{-1}(\varphi_k \mathcal{F}\psi_k).
\end{displaymath}
Here we have 
\begin{eqnarray}
\int_{\mathbb{R}^n}\frac{\partial}{\partial x} (a^tb_y^s )
\sum_{k=s+4}^{\infty}\mathcal{F}^{-1}(\varphi_k \mathcal{F}\psi_k)&=&
\int_{\mathbb{R}^n}\frac{\partial}{\partial x} (a^tb_y^s )
\mathcal{F}^{-1}\Big(\sum_{k=s+4}^{\infty}\varphi_k \mathcal{F}\psi_k\Big) \nonumber \\
&&\textrm{sine the sum is locally finite} \nonumber \\
&=&\int_{\mathbb{R}^n}\frac{\partial}{\partial x} (a^tb_y^s )\mathcal{F}\mathcal{F}^{-1}
\mathcal{F}^{-1}\Big(\sum_{k=s+4}^{\infty}\varphi_k \mathcal{F}\psi_k\Big) \nonumber \\
&=& (2\pi)^n\int_{\mathbb{R}^n} \mathcal{F}\Big( \frac{\partial }{\partial_x}(a^tb_y^s)\Big)\sum_{k=s+4}^{\infty}\varphi_k(-\;\cdot)\mathcal{F}\psi_k(-\;\cdot)\nonumber \\
&&\textrm{because} \; \frac{\partial}{\partial x}(a^tb_y^s) \in \mathcal{S} \; \textrm{and}\; \sum_{k=s+4}^{\infty}\varphi_k\mathcal{F}\psi_k) \in \mathcal{S}'\nonumber \\
&=&0.\nonumber
\end{eqnarray}
In the last step of the above calculations we used the fact that
\begin{displaymath}
 supp \; \mathcal{F}(\frac{\partial}{\partial x}(a^tb_y^s)) \subset \left\{ \xi \vert \; \vert \xi \vert \leq 5\cdot 2^s \right\}
\end{displaymath}
and
\begin{displaymath}
supp \; \sum_{k=s+4}^{\infty} \varphi_k(-  \; \cdot ) \subset \left\{\xi \vert \; 2^{s+3}\leq \vert \xi \vert\right \}
\end{displaymath}
imply that
\begin{displaymath}
supp  \;\mathcal{F}\Big(\frac{\partial}{\partial x}(a^tb_y^s)\Big) \cap supp \; \sum_{k=s+4}^{\infty}\varphi_k = \emptyset.
\end{displaymath}
This completes the proof.
\begin{flushright}
$\Box$
\end{flushright}
\newpage
Now, we can start with the estimate of $\sum_{s=0}^{\infty}\sum_{t=s-1}^{s+1}a_x^tb_y^s -a_y^tb_x^s$.
Our goal is to show that $\sum_{s=0}^{\infty}\sum_{t=s-1}^{s+1}a_x^tb_y^s -a_y^tb_x^s$ belongs to $B^0_{\mathcal{M}^{\frac{n}{2}}_1,1}$.
Making use of the above duality result, see proposition \ref{Duality}, we will first show that
\begin{displaymath}
\sum_{t=s-1}^{s+1}a_x^tb_y^s -a_y^tb_x^s \in B^0_{\mathcal{M}^{\frac{n}{2}}_1,1} \; \forall s \in \mathbb{N}
\end{displaymath} 
then we establish
\begin{displaymath}
\sum_{s=0}^{\infty} \Big \vert \Big \vert \sum_{t=s-1}^{s+1}a_x^tb_y^s -a_y^tb_x^s \Big \vert B^0_{\mathcal{M}^{\frac{n}{2}}_1,1} \Big \vert \Big\vert < \infty.
\end{displaymath}
This ensures that
\begin{displaymath}
\sum_{s=0}^{\infty}\sum_{t=s-1}^{s+1}a_x^tb_y^s -a_y^tb_x^s \in B^0_{\mathcal{M}^{\frac{n}{2}}_1,1} \subset N^0_{\frac{n}{2},1,1}.
\end{displaymath}
First of all, let us fix $t=s+j$ where $j \in \left\{-1,0,1\right\}$. \\
In order to show that $a_x^tb_y^s -a_y^tb_x^s \in B^0_{\mathcal{M}^{\frac{n}{2}}_1,1}$ it suffices to show that for all $\psi \in b^0_{L^1(H^{n-2}_\infty),\infty}$ with $\vert \vert \psi \vert b^0_{L^1(H^{n-2}_\infty),\infty} \vert \vert\leq 1$ the following inequality holds
\begin{displaymath}
\int_{\mathbb{R}^n} \psi \;d(a_x^tb_y^s -a_y^tb_x^s)  =  \int_{\mathbb{R}^n} \psi (a_x^tb_y^s -a_y^tb_x^s) \; d \lambda < \infty
\end{displaymath}
where as before $\lambda$ denotes the Lebesgue measure.\\
Moreover, in the subsequent calculations we assume that for $\psi$ we have a representation $\left\{\psi_k\right\}_{k=0}^{\infty}$, i.e. $\sum_{k=0}^{\infty}\check \varphi_k \ast \psi_k = \psi$, such that
\begin{displaymath}
\sup_{k} \vert \vert \psi_k \vert \vert _{L^1(H^{n-2}_{\infty})} \leq 2 \vert \vert \psi \vert b^0_{L^1(H^{n-2}_\infty),\infty} \vert \vert \leq 2
\end{displaymath}
and again, recall that we have density of $\mathcal{S}$ in $b^0_{L^1(H^{n-2}_\infty),\infty}$.\\
In this case we have 
\begin{eqnarray}
\int_{\mathbb{R}^n} \psi (a_x^tb_y^s-a_y^tb_x^s)
&=&\int_{\mathbb{R}^n}\psi\frac{\partial}{\partial x}\big(a^tb_y^s\big)-\psi\frac{\partial}{\partial y} \big(a^tb_x^s\big)\nonumber \\
&=&\int_{\mathbb{R}^n}\Big \lbrack\frac{\partial}{\partial x}\big(a^tb_y^s\big)\Big(\sum_{k=0}^{s+3}\mathcal{F}^{-1}(\varphi_k \mathcal{F}\psi_k) \Big ) \nonumber \\
&& -\frac{\partial}{\partial y} \big(a^tb_x^s\big)\Big(\sum_{k=0}^{s+3}\mathcal{F}^{-1}(\varphi_k \mathcal{F}\psi_k) \Big )\Big \rbrack \nonumber \\
&&\textrm{because of the same reason as in lemma \ref{PartSum}} \nonumber \\
&=&\int_{\mathbb{R}^n}\Big \lbrack-a^tb_y^s\frac{\partial}{\partial x}\Big(\sum_{k=0}^{s+3}\mathcal{F}^{-1}(\varphi_k \mathcal{F}\psi_k) \Big ) \nonumber \\
&& + a^tb_x^s\frac{\partial}{\partial y}\Big(\sum_{k=0}^{s+3}\mathcal{F}^{-1}(\varphi_k \mathcal{F}\psi_k) \Big )\Big \rbrack \nonumber \\
&&\textrm{by a simple integration by parts}\nonumber 
\end{eqnarray}
\newpage
and further
\begin{eqnarray}
\int_{\mathbb{R}^n} \psi (a_x^tb_y^s-a_y^tb_x^s)
&\leq&\int_{\mathbb{R}^n}\Big \lbrack-a^tb_y^s\Big(\sum_{k=0}^{s+3}\frac{\partial}{\partial x}\check\varphi_k \ast\psi_k \Big ) \nonumber \\
&& + a^tb_x^s\Big(\sum_{k=0}^{s+3}\frac{\partial}{\partial y}\check\varphi_k \ast\psi_k \Big )\Big \rbrack \nonumber \\
&\leq&\sum_{k=0}^{s+3}\int_{\mathbb{R}^n}\Big\lbrack -a^tb_y^s\frac{\partial}{\partial x}\check\varphi_k \ast\psi_k \nonumber \\
&& + a^tb_x^s\frac{\partial}{\partial y}\check\varphi_k \ast\psi_k\Big \rbrack \nonumber \\
&\leq& \sum_{k=0}^{s+3} \Big ( \vert \vert a^tb_y^s \vert \mathcal{M}^{\frac{n}{2}} \vert \vert \; 
\vert \vert \frac{\partial}{\partial x}\check\varphi_k \ast\psi_k \vert L^1(H^{n-2}_{\infty}) \vert \vert\nonumber \\
&& + \vert \vert a^tb_x^s \vert \mathcal{M}^{\frac{n}{2}} \vert \vert \; 
\vert \vert \frac{\partial}{\partial y}\check\varphi_k \ast\psi_k \vert L^1(H^{n-2}_{\infty}) \vert \vert \Big )\nonumber \\
&&\textrm{by proposition \ref{dualityAdams}} \nonumber\\
&\leq& \sum_{k=0}^{s+3} \Big ( \vert \vert a^tb_y^s \vert \mathcal{M}^{\frac{n}{2}}_1 \vert \vert \;  
\vert \vert \frac{\partial}{\partial x}\check\varphi_k \ast\psi_k \vert L^1(H^{n-2}_{\infty}) \vert \vert \nonumber \\
&&+ \vert \vert a^tb_x^s \vert \mathcal{M}^{\frac{n}{2}}_1 \vert \vert \; 
\vert \vert \frac{\partial}{\partial y}\check\varphi_k \ast\psi_k \vert L^1(H^{n-2}_{\infty}) \vert \vert \Big )\nonumber \\
&\leq& \sum_{k=0}^{s+3} \Big ( \vert \vert a^t \vert \mathcal{M}^{n}_2\vert \vert \; 
\vert \vert b_y^s \vert \mathcal{M}^{n}_2 \vert \vert \; 
\vert \vert \frac{\partial}{\partial x}\check\varphi_k \ast\psi_k \vert L^1(H^{n-2}_{\infty}) \vert \vert \nonumber \\
&&+ \vert \vert a^t \vert \mathcal{M}^{n}_2 \vert \vert \; 
\vert \vert b_x^s \vert \mathcal{M}^{n}_2\vert \vert \; 
\vert \vert \frac{\partial}{\partial y}\check\varphi_k \ast\psi_k \vert L^1(H^{n-2}_{\infty}) \vert \vert \Big)\nonumber \\
&&\textrm{because of H\"older's inequality with Morrey norms} \nonumber \\
&&\textrm{see also remark below}\nonumber \\
&\leq& \sum_{k=0}^{s+3} \Big ( \vert \vert a^t \vert \mathcal{M}^{n}_2\vert \vert \; 
\vert \vert b_y^s \vert \mathcal{M}^{n}_2 \vert \vert \;
2^k \vert \vert\psi\vert b^0_{L^1(H^{n-2}_{\infty}),\infty} \vert \vert \nonumber \\
&&+ \vert \vert a^t \vert \mathcal{M}^{n}_2 \vert \vert \; 
\vert \vert b_x^s \vert \mathcal{M}^{n}_2\vert \vert \;
2^k\vert \vert \psi \vert b^0_{L^1(H^{n-2}_{\infty}),\infty} \vert \vert \Big)\nonumber \\
&&\textrm{according to lemma \ref{ndimmainlemma}} \nonumber \\
&\leq& C 2^s \vert \vert a^t \vert \mathcal{M}^{n}_2\vert \vert \; 
\vert \vert b_y^s \vert \mathcal{M}^{n}_2 \vert \vert 
+ C 2^s \vert \vert a^t \vert \mathcal{M}^{n}_2 \vert \vert \; 
\vert \vert b_x^s \vert \mathcal{M}^{n}_2\vert \vert \nonumber \\
&<& \infty \nonumber \\
&&\textrm{due to our assumptions.} \nonumber
\end{eqnarray}
Thus we have seen that for all $s \in \mathbb{N}$
\begin{displaymath}
a_x^tb_y^s-a_y^tb_x^s \in (b^0_{L^1(H^{n-2}_\infty),\infty})^* = B^0_{\mathcal{M}^{\frac{n}{2}}_1,1} \subset N^0_{\frac{n}{2},1,1}.
\end{displaymath}
Next, we study
\begin{displaymath}
\sum_{s=0}^{\infty}\Big\vert \Big\vert \sum_{t=s-1}^{s+1}a_x^tb_y^s -a_y^tb_x^s \Big\vert B^0_{\mathcal{M}^{\frac{n}{2}}_1,1}\Big\vert \Big\vert.
\end{displaymath}
What concerns this latter quantity, we will assume for the sake of simplicity that $t=s$. Then we can estimate
\begin{eqnarray}
\sum_{s=0}^{\infty}\vert \vert a_x^sb_y^s -a_y^sb_x^s \vert B^0_{\mathcal{M}^{\frac{n}{2}}_1,1}\vert \vert &=&
\vert \vert a_x^0b_y^0 -a_y^0b_x^0 \vert B^0_{\mathcal{M}^{\frac{n}{2}}_1,1}\vert \vert
+ \sum_{s=1}^{\infty}\vert \vert a_x^sb_y^s -a_y^sb_x^s \vert B^0_{\mathcal{M}^{\frac{n}{2}}_1,1}\vert \vert \nonumber \\
&\leq& C \vert \vert a^0 \vert \mathcal{M}^{n}_2\vert \vert \; 
\vert \vert b_y^0 \vert \mathcal{M}^{n}_2 \vert \vert 
+ C  \vert \vert a^0 \vert \mathcal{M}^{n}_2 \vert \vert \; 
\vert \vert b_x^0 \vert \mathcal{M}^{n}_2\vert \vert \nonumber \\
&&+ C \sum_{s=1}^{\infty} 2^s \vert \vert a^s \vert \mathcal{M}^{n}_2\vert \vert \; 
\vert \vert b_y^s \vert \mathcal{M}^{n}_2 \vert \vert \nonumber \\
&& + C \sum_{s=1}^{\infty} 2^s \vert \vert a^s \vert \mathcal{M}^{n}_2 \vert \vert \; 
\vert \vert b_x^s \vert \mathcal{M}^{n}_2\vert \vert \nonumber \\
&\leq& C \vert \vert a^0 \vert \mathcal{M}^{n}_2\vert \vert \; 
\vert \vert b_y^0 \vert \mathcal{M}^{n}_2 \vert \vert 
+ C  \vert \vert a^0 \vert \mathcal{M}^{n}_2 \vert \vert \; 
\vert \vert b_x^0 \vert \mathcal{M}^{n}_2\vert \vert \nonumber \\
&&+ C \sum_{s=1}^{\infty} \vert \vert a^s_x \vert \mathcal{M}^{n}_2\vert \vert \; 
\vert \vert b_y^s \vert \mathcal{M}^{n}_2 \vert \vert \nonumber \\
&& + C \sum_{s=1}^{\infty} \vert \vert a^s_y \vert \mathcal{M}^{n}_2 \vert \vert \; 
\vert \vert b_x^s \vert \mathcal{M}^{n}_2\vert \vert \nonumber \\
&&\textrm{similar to $2^{ms}\vert \vert g \vert \vert_p \simeq \vert \vert \nabla^m g \vert \vert_p$ (under appropriate assumptions)} \nonumber \\
&&\textrm{cf. also theorem 2.9 in \cite{KY}}\nonumber \\
&\leq& C \vert \vert a^0 \vert \mathcal{M}^{n}_2\vert \vert \; 
\vert \vert b_y^0 \vert \mathcal{M}^{n}_2 \vert \vert 
+ C  \vert \vert a^0 \vert \mathcal{M}^{n}_2 \vert \vert \; 
\vert \vert b_x^0 \vert \mathcal{M}^{n}_2\vert \vert \nonumber \\
&&+ C \Big(\sum_{s=1}^{\infty} \vert \vert a^s_x \vert \mathcal{M}^{n}_2\vert \vert^2 \Big)^{\frac{1}{2}}  
\Big( \sum_{s=1}^{\infty} \vert \vert b_y^s \vert \mathcal{M}^{n}_2 \vert \vert^2 \Big)^{\frac{1}{2}} \nonumber \\
&& + C \Big(\sum_{s=1}^{\infty} \vert \vert a^s_y \vert \mathcal{M}^{n}_2 \vert \vert^2 \Big)^{\frac{1}{2}}
\Big(\sum_{s=1}^{\infty}\vert \vert b_x^s \vert \mathcal{M}^{n}_2\vert \vert ^2\Big)^{\frac{1}{2}}\nonumber \\
&&\textrm{by H\"older's inequality} \nonumber \\
&<&\infty \nonumber \\
&&\textrm{thanks to our hypothesis.} \nonumber
\end{eqnarray}
All together we have seen that 
\begin{displaymath}
\sum_{s=0}^{\infty} a_x^sb_y^s -a_y^sb_x^s \in B^0_{\mathcal{M}^{\frac{n}{2}}_1,1}\subset N^0_{\frac{n}{2},1,1}.
\end{displaymath}

Now, since the above estimate is independent of the choice of $j$ we immediately conclude that
\begin{displaymath}
\sum_{s=0}^{\infty}\sum_{t=s-1}^{s+1}a_x^tb_y^s -a_y^tb_x^s \in N^0_{\frac{n}{2},1,1}
\end{displaymath}

Now, as we know that $\sum_{s=0}^{\infty}\sum_{t=s-1}^{s+1}a_x^tb_y^s -a_y^tb_x^s \in B^0_{\mathcal{M}^{\frac{n}{2}}_1,1} \subset N^0_{\frac{n}{2},1,1}$ we apply the embedding result of Kozono/Yamazaki, theorem 2.5 in \cite{KY}, and find that
\begin{displaymath}
\sum_{s=0}^{\infty}\sum_{t=s-1}^{s+1}a_x^tb_y^s -a_y^tb_x^s \in B^{-2}_{\infty,1}.
\end{displaymath}

\begin{Remark}
\textnormal{
Assume that $f,g\in \mathcal{M}^n_2$. Then we have for all $0<r$ and for all $x \in \mathbb{R}^n$
\begin{eqnarray}
\vert \vert fg \vert \vert_{L^1(B_r(x))} &\leq& \vert \vert f \vert \vert _{L^2(B_r(x))} \vert \vert g \vert \vert _{L^2(B_r(x))} \nonumber \\
&\leq& C_1 r^{\frac{n}{2}- 1} C_2 r^{\frac{n}{2}- 1}\nonumber \\
&=& C r^{n-2}.\nonumber 
\end{eqnarray}
According to the definition, this shows that $fg \in \mathcal{M}^{\frac{n}{2}}_1$.
}
\end{Remark}

\textbf{\textit{Regularity}}\\

We rewrite our equation $\Delta u=f$ as $\Delta u = f^0 +\sum_{k\geq 1} f^k$. 

And the solution $u$ can be written as
\begin{eqnarray}
u&=& \Delta ^{-1} f^0 + \Delta ^{-1} (\sum_{k\geq 1} f^k )\nonumber \\
&=:& u_1 +u_2. \nonumber
\end{eqnarray}

Our strategy is to show that $u_1$ as well as $u_2$ is continuous and bounded. \\

What concerns $u_1$, observe that due to the Paley-Wiener theorem $f^0$ is analytic, so in particular continuous. This implies immediately - by classical results (see e.g. \cite{GT}) - that $u_1$ is continuous. \\
On one hand we have that $f^0 \in B^s_{\frac{n}{2},2}$ for all $s\in \mathbb{R}$ (since $\nabla a, \nabla b \in B^0_{\mathcal{M}^n_2,2} \subset \mathcal{M}^n_2 \subset L^n$) on the other hand we know that $f^0 \in B^s_{\infty,1}$ for all $s\in \mathbb{R}$ because $f \in B^{-2}_{\infty,1}$. From that we can deduce by standard elliptic estimates (see also \cite{RS}) and the embedding result of Sickel and Triebel \cite{ST} that $u_1$ is not only continuous but also bounded!

Next, we will show that $u_2$ is bounded and continuous. In order to reach this goal, we show that $u_2 \in B^0_{\infty,1}$: We find the following estimates
\begin{eqnarray}
\vert \vert u_2 \vert B^0_{\infty,1} \vert \vert &=&  \sum_{s=0}^{\infty}\vert \vert  u_2^s \vert \vert_{\infty} \nonumber \\
&=&  \sum_{s=0}^{\infty}2^{-2s} 2^{2s}\vert \vert  u_2^s \vert \vert_{\infty} \nonumber \\
&=&C \sum_{s=0}^{\infty} 2^{-2s}\vert \vert(\Delta u_2)^s \vert \vert _\infty \nonumber 
\end{eqnarray}

This last passage holds thanks to the fact that 
\begin{displaymath}\label{tao}
2^{ms}\vert \vert g \vert \vert_p \simeq \vert \vert \nabla^m g \vert \vert_p
\end{displaymath} 
if the Fourier transform of $g$ is supported on an annulus with radii comparable to $2^s$(see \cite{Tao} for instance). 

For $s=0$ we observe
\begin{displaymath}
\mathcal{F}( -\Delta u_2 ) = \mathcal{F}( \sum_{k\geq 1} f^k)
\end{displaymath}
which implies
\begin{displaymath}
supp (\mathcal{F}(u_2)) \subset  (B_1(0))^c
\end{displaymath}
because of the fact that 
\begin{displaymath}
supp ( \mathcal{F}(\sum_{k\geq 1} f^k)) \subset (B_1(0))^c.
\end{displaymath}
So in this case too, we can apply the above mentioned fact in order to conclude that also for $s=0$ we have
\begin{displaymath}
\vert \vert u_2^0 \vert \vert_{\infty} \leq C \vert \vert (\Delta u_2)^0 \vert \vert_{\infty}.
\end{displaymath}

Back to our estimate, we continue
\begin{eqnarray}
\vert \vert u_2 \vert B^0_{\infty,1} \vert \vert &\leq& C \sum_{s=0}^{\infty}2^{-2s}\vert \vert (\Delta u_2)^s \vert \vert _{\infty} \nonumber \\
&=& C \sum_{s=0}^{\infty}2^{-2s}\vert \vert (\sum_{k\geq1}f^k)^s \vert \vert _{\infty} \nonumber \\
&=& C\sum_{s=0}^{\infty}2^{-2s}\vert \vert \mathcal{F}^{-1}(\sum_{k=s-1}^{s+1}\varphi_s \varphi_k \hat f) \vert \vert _{\infty} \nonumber \\
&\leq& \sum_{s=0}^{\infty}2^{-2s}\vert \vert f^s \vert \vert _{\infty} \nonumber \\
&&\textrm{thanks to a Fourier multiplier result} \nonumber \\
&&\textrm{for further details we refer to \cite{T}} \nonumber \\
&=& C \vert \vert f \vert B^{-2}_{\infty,1} \vert \vert \nonumber \\
&<& \infty \; \textrm{according to our assumptions.} \nonumber 
\end{eqnarray}

This shows that $u_2$ belongs to $B^0_{\infty,1}(\mathbb{R}^n)$. 

Alternatively one could make use of the lifting property, see \cite{RS}, chapter 2.6, to show that $u_2 \in C$. (Recall that $C$ denotes the space of all uniformly continuous functions on $\mathbb{R}^n$.)
The last ingredient is the embedding result due to Sickel/Triebel (see \cite{ST}). 

This leads immediately to the assertion we claimed because $u$ as a sum of two bounded continuous functions is again continuous and bounded.
\begin{flushright}
$\Box$
\end{flushright}

\subsection{Proof of theorem \ref{regndim} ii)}

In a first step we show that $a_xb_y-a_yb_x \in B^{-1}_{\mathcal{M}^n_2,1}$: \\
From the proof of theorem \ref{regndim} we know that
\begin{displaymath}
\sum_{k=0}^{\infty} \sum_{s=k-1}^{k+1} a_x^k b_y^s -a_y^k b_x^s \in B^0_{\mathcal{M}^{\frac{n}{2}}_1,1} \subset B^{-1}_{\mathcal{M}^n_2,1}.
\end{displaymath}
Next, we observe that
\begin{eqnarray}
\vert \vert \pi_3(a_x,b_x) \vert B^{-1}_{\mathcal{M}^n_2,1}\vert \vert &\leq& C\sum_{s=0}^{\infty} 2^{-s} \Big\vert \Big \vert \sum_{k=0}^{s-2} a_x^sb_y^k \Big\vert  \Big\vert_{\mathcal{M}^n_2} \nonumber \\
&&\textrm{by a simple modification of lemma 3.16 in \cite{Maz2}} \nonumber \\
&\leq& C \sum_{s=0}^{\infty} 2^{-s} \vert  \vert a_x^s \vert \vert _{\mathcal{M}^n_2} \Big \vert \Big \vert \sum_{k=0}^{s-2} b_y^k \Big\vert  \Big\vert_{\infty} \nonumber \\
&\leq& C \Big( \sum_{s=0}^{\infty} \vert \vert a_x^s \vert \vert_{\mathcal{M}^n_2}^2 \Big )^{\frac{1}{2}} \Big ( \sum_{s=0}^{\infty}2^{-2s} \Big \vert \Big \vert \sum_{k=0}^{s-2} b_y^k \Big \vert \Big \vert_{\infty}^2 \Big )^{\frac{1}{2}} \nonumber \\
&\leq& C \vert \vert a_x \vert B^0_{\mathcal{M}^n_2,2} \vert \vert \Big(\sum_{s=0}^{\infty} 2^{-2s} \Big \vert \Big \vert \sum_{k=0}^sb_y^k \Big \vert \Big \vert_{\infty}^2 \Big)^{\frac{1}{2}}\nonumber \\
&\leq&C\vert \vert a_x \vert B^0_{\mathcal{M}^n_2,2} \vert \vert \; \vert \vert b_y \vert B^{-1}_{\mathcal{M}^n_2,2}\vert \vert \nonumber \\
&&\textrm{according to lemma 4.4.2 of \cite{RS}} \nonumber \\
&\leq& \vert \vert a_x \vert B^0_{\mathcal{M}^n_2,2} \vert \vert \; \vert \vert b_y \vert B^0_{\mathcal{M}^n_2,2}\vert \vert. \nonumber 
\end{eqnarray}
Now, since
\begin{displaymath}
\partial _{x_i}u= \mathcal{F}^{-1}\Big ( i \frac{\xi_i}{\vert \xi \vert^2}\mathcal{F}(\Delta u)\Big)
\end{displaymath}
we note first, that due to the facts that $\Delta u \in F^0_{1,2} \subset L^1$ and $r^{-1} \in L^{\frac{n}{n-1}}$ for $n\geq 3$,
\begin{displaymath}
(\nabla u)^0 \in L^n \subset \mathcal{M}^n_2
\end{displaymath}
which implies that $(\nabla u)^0 \in B^0_{\mathcal{M}^n_2,2}$.\\
Second, for $s \geq 1$ we have
\begin{displaymath}
\vert \vert (\nabla u)^s \vert \vert_{\mathcal{M}^n_2} \leq C 2^{-s} \vert \vert (\Delta u)^s \vert \vert_{\mathcal{M}^n_2} 
\end{displaymath}
which leads to the conclusion - remember the first step! - that $\sum_{s\geq 1}(\nabla u)^s \in B^0_{\mathcal{M}^n_2,1}$. \\
Alternatively one could observe that 
\begin{displaymath}
\Big \vert \partial^{\vert \alpha \vert} \Big(\frac{\xi_i}{\vert \xi \vert^2}\Big ) \Big \vert \leq C \vert \xi \vert ^{-1-\vert \alpha \vert}
\end{displaymath}
information, which together with theorem 2.9 in \cite{KY} leads to the same conclusion as above, namely that
\begin{displaymath}
\nabla u \in B^0_{\mathcal{M}^n_2,1}.
\end{displaymath}

These estimates complete the proof.
\begin{flushright}
$\Box$
\end{flushright}

\subsection{Proof of theorem \ref{regndim} iii)}

This proof is very similar to the one of theorem \ref{regndim} ii).\\
In stead of the observation $\Big \vert \partial^{\vert \alpha \vert} \Big(\frac{\xi_i}{\vert \xi \vert^2}\Big )\Big \vert \leq C \vert \xi \vert ^{-1-\vert \alpha \vert}$ here we use theorem 2.9 of \cite{KY} together with the fact that
\begin{displaymath}
\Big \vert \partial^{\vert \alpha \vert} \Big(\frac{\xi_i \xi_j}{\vert \xi \vert^2}\Big ) \Big \vert\leq C \vert \xi \vert ^{-\vert \alpha \vert}.
\end{displaymath}
\begin{flushright}
$\Box$
\end{flushright}

\subsection{Proof of theorem \ref{gauge}}

\begin{Lemma}\label{lemmagauge}
There exist constants $\varepsilon(m)>0$ and $C(m) >0$ such that for every $\Omega \in B^0_{\mathcal{M}^n_2,2}(B^n_1(0), so(m)\otimes \Lambda^1\mathbb{R}^n)$ which satisfies
\begin{displaymath} 
\vert \vert \Omega \vert B^0_{\mathcal{M}^n_2,2} \vert \vert \leq \varepsilon(m)
\end{displaymath}
there exist $\xi \in B^1_{\mathcal{M}^n_2,2}(B^n_1(0),so(m)\otimes\Lambda^{n-2}\mathbb{R}^n)$ and $P \in B^1_{\mathcal{M}^n_2,2}(B^n_1(0),SO(m))$ such that
\begin{itemize}
\item[i)]
\begin{displaymath}
*d\xi= P^{-1}dP+P^{-1}\Omega P \; \textrm{in} \; B^n_1(0)
\end{displaymath}
\item[ii)]
\begin{displaymath}
\xi=0 \; \textrm{on} \; \partial B^n_1(0)
\end{displaymath}
\item[iii)]
\begin{displaymath}
\vert \vert \xi \vert B^1_{\mathcal{M}^n_2,2} \vert \vert + \vert \vert P \vert B^1_{\mathcal{M}^n_2,2} \vert \vert\leq C(m)\vert \vert \Omega \vert B^0_{\mathcal{M}^n_2,2} \vert \vert.
\end{displaymath}
\end{itemize}
\end{Lemma}

The proof of this lemma is a straightforward adaptation of the corresponding assertion in \cite{RivSt}.\\

Now, let $\varepsilon(m), P$ and $\xi$ be as in lemma \ref{lemmagauge}. Note that since $P\in SO(m)$ we have also $P^{-1} \in B^1_{\mathcal{M}^n_2,2}$. Our goal is to find $A$ and $B$ such that
\begin{equation}\label{glAB}
dA-A\Omega=-d^*B.
\end{equation} 
If we set $\tilde{A}:=AP$ then, according to equation (\ref{glAB}) it has to satisfy
\begin{displaymath}
d\tilde{A}+(d^*B)P=\tilde{A}+d\xi.
\end{displaymath}
As a intermediate step we will first study the following problem
\begin{eqnarray}\nonumber
\left\{ \begin{array}{ll}
\Delta \hat{A} &= d\hat{A} \cdot *d\xi -d^* B \cdot \nabla P \; \textrm{in} B^n_1(0)\\
d(d^* B) &= d\hat{A}\wedge dP^{-1} -d*(\hat{A}d\xi P^{-1})-d*(d\xi P^{-1}) \\
\frac{\partial \hat{A}}{\partial \nu}, &= 0 \; \textrm{and} \; B=0 \; \textrm{on} \; \partial B^n_1(0) \\
\int_{B^n_1(0)}\hat{A} &= id_m.
\end{array} \right.
\end{eqnarray}
For this system we have the a-priori-estimates (recall theorem \ref{regndim} with its proof, lemma \ref{stab} and the fact that we are working on a bounded domain) 
\begin{eqnarray}
\vert \vert  \hat{A} \vert B^1_{\mathcal{M}^n_2,2} \vert \vert + \vert \vert \hat{A} \vert \vert_{\infty} 
&\leq& C \vert \vert \xi \vert B^1_{\mathcal{M}^n_2,2} \vert \vert \; \vert \vert \hat{A} \vert B^1_{\mathcal{M}^n_2,2} \vert \vert \nonumber \\
&&+ C \vert \vert P \vert B^1_{\mathcal{M}^n_2,2} \vert \vert \; \vert \vert B \vert B^1_{\mathcal{M}^n_2,2} \vert \vert \nonumber
\end{eqnarray}
and
\begin{eqnarray}
\vert \vert B \vert B^1_{\mathcal{M}^n_2,2} \vert \vert 
&\leq& C \vert \vert P^{-1} \vert B^1_{\mathcal{M}^n_2,2} \vert \vert \; \vert \vert \hat{A} \vert B^1_{\mathcal{M}^n_2,2} \vert \vert 
+ C \vert \vert \xi \vert B^1_{\mathcal{M}^n_2,2} \vert \vert \; \vert \vert \hat{A} \vert \vert_{\infty} \nonumber \\
&&+C \vert \vert \xi \vert B^1{\mathcal{M}^n_2,2} \vert \vert. \nonumber
\end{eqnarray}
Since the used norms of $\xi$ and $P$ - as well as of $P^{-1}$ - can be bounded in terms of $C\vert \vert \Omega \vert B^0_{\mathcal{M}^n_2,2} \vert \vert$ the above estimates together with standard fixpoint theory guarantee the existence of $\hat{A}$ and $B$ such that they solve the above system and in addition satisfy
\begin{equation}\label{estimate}
\vert \vert \hat{A} \vert B^1_{\mathcal{M}^n_2,2} \vert \vert + \vert \vert \hat{A} \vert \vert_{\infty} + \vert \vert B \vert B^1_{\mathcal{M}^m_2,2} \vert \vert
\leq C \vert \vert \Omega \vert B^0_{\mathcal{M}^m_2,2} \vert \vert.
\end{equation}
Next, similar to the proof of corollary \ref{ndimRiv} we decompose for some $D$
\begin{displaymath}
d\hat{A}-\hat{A}*d\xi +d^*BP=d^*D.
\end{displaymath}
Then we set $\tilde{A}:=\hat{A}+id_m$, which satisfies for some $n-2$-form $F$
\begin{displaymath}
d\tilde{A}-\tilde{A}*d\xi +d^*B P=d^*D -*d\xi =: *dF.
\end{displaymath}
It is not difficult to show that $*d(*dFP^{-1})=0$ together with $F=0$ on $\partial B^n_1(0)$ imply that $F\equiv 0$ (see also a similar assertion in \cite{Riv} and remember that on bounded domains $B^0_{\mathcal{M}^n_2,2} \subset L^2$). \\
From this we conclude that in fact $\tilde{A}$ satisfies the desired equation. If wet finally set $A:=\tilde{A}P^{-1}$ and let $B$ as given in the above system we get that in fact these $A$ and $B$ solve the required relation (\ref{glAB}). \\
So far, we have proved parts ii) and iii) of theorem \ref{gauge} (recall also estimate (\ref{estimate})).\\
Moreover, the invertibility of $A$ follows immediately from its construction, likewise the estimates for $\nabla A$ and $\nabla A^{-1}$.\\
Last but not least, the relation $A=\hat{A}P^{-1} +id_mP^{-1}$ implies that
\begin{displaymath}
\vert \vert dist(A,SO(m))\vert \vert_{\infty}\leq C \vert \vert \hat{A} \vert \vert_{\infty} \leq C\vert \vert \Omega \vert B^0_{\mathcal{M}^n_2,2} \vert \vert.
\end{displaymath}
This completes the proof of theorem \ref{gauge}.
\begin{flushright}
$\Box$
\end{flushright}

\subsection{Proof of corollary \ref{ndimRiv}}

The first part of the corollary is a straightforward calculation.\\
Let $A$ and $B$ be as in theorem \ref{gauge}.\\
Then we have
\begin{displaymath}
\left\{ \begin{array}{ll}
*d*(Adu) = -d^*B\cdot \nabla u \\
d(Adu) =  dA\wedge du.
\end{array} \right.
\end{displaymath}
These equations together with a classical Hodge decomposition for $Adu$
\begin{displaymath}
Adu=d^*E+dD \; \textrm{with} E,D \in W^{1,2}
\end{displaymath}
lead to the following equations
\begin{displaymath}
\left\{ \begin{array}{ll}
\Delta D = -d^*B\cdot \nabla u \\
\Delta E =  dA\wedge du.
\end{array} \right.
\end{displaymath}
Since the right hand sides are made of Jacobians we conclude that $D,E \in B^0_{\infty,1}$. Next, we observe that
\begin{displaymath}
du=A^{-1}(d^*E+dD) \in B^0_{\mathcal{M}^n_2,1} \subset B^{-1}_{\infty,1}. 
\end{displaymath}
This holds because $A^{-1}\in B^1_{\mathcal{M}^n_2,2} \cap L^{\infty}$ (see also theorem \ref{gauge}) and $dD, d^*E \in B^0_{\mathcal{M}^n_2,1}$ (see also theorem \ref{regndim} ii)).The proof of the above fact is the same as the proof of the assertion of lemma \ref{stab}. In a last step we note that (recall the reasons why theorem \ref{regndim} hold) thanks to the information we have so far 
\begin{displaymath}
u \in B^0_{\infty,1} \subset C
\end{displaymath}
which completes the proof.
\begin{flushright}
$\Box$
\end{flushright}

\subsection{Proof of lemma \ref{inMorrey}}

We start with the following observation. \\
Let $x_0 \in \mathbb{R}^n$ and $r > 0$ and recall that $1< q \leq 2$ and $r \leq q$. Then for $f \in B^0_{\mathcal{M}^p_q,r}$ we have
\begin{eqnarray}
\Big( \int_{B_r(x_0)} \big( \sum_{s=0}^{\infty} \vert f^s \vert ^2\big)^{\frac{q}{2}} \Big )^{\frac{1}{q}} 
&\leq& \Big( \int_{B_r(x_0)} \sum_{s=0}^{\infty} \vert f^s \vert ^q \Big )^{\frac{1}{q}} \nonumber \\
&\leq& \Big( \sum_{s=0}^{\infty}\int_{B_r(x_0)}  \vert f^s \vert ^q \Big )^{\frac{1}{q}} \nonumber \\
&\leq& \Big( \sum_{s=0}^{\infty}\vert \vert  f^s \vert \vert^q_{L^q(B_r(x_0))} \Big )^{\frac{1}{q}} \nonumber \\
&\leq& \Big( \sum_{s=0}^{\infty}\vert \vert  f^s \vert \vert^q_{\mathcal{M}^p_q}(r^{\frac{n}{q}-\frac{n}{p}})^q \Big )^{\frac{1}{q}} \nonumber \\
&\leq& \Big( (r^{\frac{n}{q}-\frac{n}{p}})^q\sum_{s=0}^{\infty} \vert \vert  f^s \vert \vert^q_{\mathcal{M}^p_q} \Big )^{\frac{1}{q}} \nonumber \\
&=& r^{\frac{n}{q}-\frac{n}{p}}\Big( \sum_{s=0}^{\infty} \vert \vert  f^s \vert \vert^q_{\mathcal{M}^p_q} \Big )^{\frac{1}{q}} \nonumber \\
&=& r^{\frac{n}{q}-\frac{n}{p}} \vert \vert f \vert B^0_{\mathcal{M}^p_q,q} \vert \vert \nonumber \\
&\leq& C r^{\frac{n}{q}-\frac{n}{p}} \vert \vert f \vert B^0_{\mathcal{M}^p_q,r} \vert \vert. \nonumber 
\end{eqnarray} 
From the last inequality we have that for all $r > 0$ and for all $x_0 \in \mathbb{R}^n$
\begin{displaymath}
r^{\frac{n}{p}-\frac{n}{q}} \Big\vert\Big \vert \big( \sum_{s=0}^{\infty} \vert f^s \vert ^2\big)^{\frac{q}{2}} \Big\vert \Big\vert_{L^q(B_r(x_0))} \leq 
C \vert \vert f \vert B^0_{\mathcal{M}^p_q,r} \vert \vert.
\end{displaymath}
This last estimate together \cite{Maz1}, proposition 4.1, implies that $f \in \mathcal{M}^p_q$. \\
The assertion in the case $f \in N^0_{p,q,r}$ is the same.\\
\begin{flushright}
$\Box$
\end{flushright}

\subsection{Proof of lemma \ref{aequiv}}

\begin {itemize}
\item[i)]
In a first step we will show that if $f\in B^1_{\mathcal{M}^p_q,r}$ there exist a constant $C$ - independent of $f$ - such that
\begin{displaymath}
\vert \vert f \vert B^0_{\mathcal{M}^p_q,r} \vert \vert + \vert \vert \nabla f \vert B^0_{\mathcal{M}^p_q,r} \vert \vert 
\leq C \vert \vert f \vert B^1_{\mathcal{M}^p_q,r} \vert \vert.
\end{displaymath}
Obviously, we have that
\begin{displaymath}
\vert \vert f \vert B^0_{\mathcal{M}^p_q,r} \vert \vert \leq \vert \vert f \vert B^1_{\mathcal{M}^p_q,r} \vert \vert.
\end{displaymath}
Moreover, we observe that
\begin{eqnarray}
\vert \vert \nabla f \vert B^0_{\mathcal{M}^p_q,r} \vert \vert &=& \Big ( \sum_{j=0}^{\infty} \vert \vert ( \nabla f )^j \vert \vert_{\mathcal{M}^q_p}^r \Big)^{\frac{1}{r}} \nonumber \\
&\leq& \Big ( \sum_{j=1}^{\infty} \vert \vert ( \nabla f )^j \vert \vert_{\mathcal{M}^q_p}^r \Big)^{\frac{1}{r}} + \vert \vert (\nabla f)^0 \vert \vert_{\mathcal{M}^p_q}  \nonumber \\
&\leq& C \Big ( \sum_{j=1}^{\infty} 2^{jr}\vert \vert f^j \vert \vert_{\mathcal{M}^q_p}^r \Big)^{\frac{1}{r}} + C \vert \vert f \vert \vert_{\mathcal{M}^p_q} \nonumber \\
&&\textrm{where for the first addend we used an estimate similar to (\ref{tao})} \nonumber \\
&&\; \; \; \textrm{with the necessary adaptations to our situation} \nonumber \\
&&\; \; \; \textrm{see also \cite{KY}}\nonumber \\
&&\textrm{and for the second addend we used \cite{KY}, lemma 1.8} \nonumber \\
&&\; \; \; \textrm{and the observation $\mathcal{F}^{-1}(\xi \varphi_0 \hat{f})=\mathcal{F}^{-1}(\xi \varphi_0)\ast f$.} \nonumber \\
&\leq& C \vert \vert f \vert B^1_{\mathcal{M}^p_q,r} \vert \vert + C \vert \vert f \vert B^0_{\mathcal{M}^p_q,r}\vert \vert \nonumber \\
&&\textrm{because of lemma \ref{inMorrey}} \nonumber \\
&\leq& C \vert \vert f \vert B^1_{\mathcal{M}^p_q,r} \vert \vert + C \vert \vert f \vert B^1_{\mathcal{M}^p_q,r}\vert \vert \nonumber \\
&\leq& \vert \vert f \vert B^1_{\mathcal{M}^p_q,r} \vert \vert\nonumber 
\end{eqnarray}
as desired.
\item[ii)]
Now, we assume that $f$ satisfies
\begin{displaymath}
\vert \vert f \vert B^0_{\mathcal{M}^p_q,r} \vert \vert + \vert \vert \nabla f \vert B^0_{\mathcal{M}^p_q,r} \vert \vert < \infty.
\end{displaymath}
We have to show that this last quantity controls
\begin{displaymath}
\vert \vert f \vert B^1_{\mathcal{M}^p_q,r} \vert \vert.
\end{displaymath}
In fact, we calculate
\begin{eqnarray}
\vert \vert f \vert B^1_{\mathcal{M}^p_q,r} \vert \vert &=& 
\Big ( \sum_{j=0}^{\infty} 2^{jr}\vert \vert f ^j \vert \vert_{\mathcal{M}^q_p}^r \Big)^{\frac{1}{r}} \nonumber \\
&\leq& C \vert \vert f^0 \vert \vert_{\mathcal{M}^p_q} + C\Big ( \sum_{j=1}^{\infty} 2^{jr}\vert \vert f^j \vert \vert_{\mathcal{M}^q_p}^r \Big)^{\frac{1}{r}} \nonumber \\
&\leq& C \vert \vert f^0 \vert B^0_{\mathcal{M}^p_q,r} + C \vert \vert \nabla f \vert B^0_{\mathcal{M}^p_q,r} \vert \vert \nonumber \\
&&\textrm{again by an adaption of estimate (\ref{tao})} \nonumber \\
&\leq& C (\vert \vert f^0 \vert B^0_{\mathcal{M}^p_q,r} + \vert \vert \nabla f \vert B^0_{\mathcal{M}^p_q,r} \vert \vert ). \nonumber 
\end{eqnarray}
\begin{flushright}
$\Box$
\end{flushright}
\end{itemize}

\subsection{Proof of lemma \ref{aufKomp}}

According to lemma \ref{aequiv} it is enough to show that $f\in B^0_{\mathcal{M}^p_q,r}$. First of all, we observe that
\begin{eqnarray}
\Big \vert \Big \vert \sum_{j=1}^{\infty} f^j \Big \vert B^0_{\mathcal{M}^p_q,r} \Big\vert \Big \vert &\leq& 
\Big \vert \Big \vert \sum_{j=1}^{\infty} f^j \Big \vert B^1_{\mathcal{M}^p_q,r} \Big \vert \Big \vert \nonumber \\
&\leq& C \Big( \sum_{j=0}^{\infty}  2^{jr} \vert \vert f^j \vert \vert_{\mathcal{M}^p_q}^r \Big)^{\frac{1}{r}} \nonumber \\
&\leq& C \Big( \sum_{j=0}^{\infty} \vert \vert (\nabla f) ^j \vert \vert_{\mathcal{M}^p_q}^r \Big)^{\frac{1}{r}} \nonumber \\
&\leq& \vert \vert \nabla f\vert B^0_{\mathcal{M}^p_q,r} \vert \vert. \nonumber
\end{eqnarray}
Now, it remains to estimate $\vert \vert f^0 \vert \vert_{\mathcal{M}^p_q}$:\\
It holds 
\begin{displaymath}
f^0=\mathcal{F}^{-1}\big (\sum_{i=1}^n\frac{\xi_i}{\vert \xi \vert^2}\xi_i \hat{f} \varphi_0 \big).
\end{displaymath}
Next, due to lemma \ref{inMorrey} and its corollary we know that $f \in L^q$ and in particular - since $f$ has compact support $f\in L^1$ so $\xi_i \hat{f} \in L^{\infty}$ for all $i$. Moreover, thanks to our assumptions
\begin{displaymath}
\varphi_0 \frac{1}{\vert \xi \vert} \in L^{\frac{p}{p-1}} \; \; \; \textrm{where} \;\; \; \frac{p}{p-1 } \in \lbrack 1,2 \rbrack.
\end{displaymath}
So, for all possible $i$
\begin{displaymath}
\varphi_0 \frac{\xi_i}{\vert \xi \vert^2}\xi_i \hat{f} \in L^{\frac{p}{p-1}}.
\end{displaymath}
From this we conclude that
\begin{displaymath}
f^0 \in L^p \subset \mathcal{M}^p_q,
\end{displaymath}
and finally
\begin{eqnarray}
\vert \vert f^0\vert B^0_{\mathcal{M}^p_q,r} \vert \vert &\leq& \vert \vert f^0 \vert \vert_{\mathcal{M}^p_q} +\vert \vert f^1\vert \vert_{\mathcal{M}^p_q} \nonumber \\
&\leq&  \vert \vert f^0 \vert \vert_{L^p} + C\Big\vert \Big \vert \sum_{j=1}^{\infty} f^j\vert B^0_{\mathcal{M}^p_q,r} \Big\vert\Big \vert\nonumber \\
&\leq& C \vert \vert \nabla f \vert B^0_{\mathcal{M}^p_q,r} \vert \vert +C \Big\vert \Big \vert \sum_{j=1}^{\infty} f^j\vert B^0_{\mathcal{M}^p_q,r} \Big \vert \Big \vert \nonumber \\
&\leq& C \vert \vert \nabla f \vert B^0_{\mathcal{M}^p_q,r} \vert \vert.\nonumber
\end{eqnarray}
\begin{flushright}
$\Box$
\end{flushright}

\subsection{Proof of lemma \ref{density}}

\textit{Density of $O_M$ in $N^s_{p,q,r}$ respectively in $B^s_{\mathcal{M}^p_q,r}$}\\
The idea is to approximate $f \in N^s_{p,q,r}$ by $f_n:= \sum_{k=0}^n f^k$.

From the definition of the spaces $N^s_{p,q,r}$ we immediately deduce that there exists $N \in \mathbb{N}$ such that 
\begin{displaymath}
\Big ( \sum_{j=N+1}^{\infty} 2^{sjr} \vert \vert f^j \vert \vert_{M^p_q}^r \Big)^{\frac{1}{r}} < \varepsilon.
\end{displaymath}
What concerns the first contributions, i.e. $f^0$-$f^N$, we know that 
\begin{displaymath}
\sum_{j=0}^{N} f^j=:f_N \in O_M.
\end{displaymath}
So, 
\begin{displaymath}
\vert \vert f- f_N \vert N^s_{p,q,r} \vert \vert \leq C \Big ( \sum_{j=N+1}^{\infty} 2^{sjr} \vert \vert f^j \vert \vert_{M^p_q}^r \Big)^{\frac{1}{r}} < C \varepsilon
\end{displaymath}
where $C$ does not depend on $f$. This shows that $f_N$ approximates $f$ in the desired way.\\
The proof in the case $B^s_{\mathcal{M}^p_q,r}$ is the same - with the necessary modifications of course.

\textit{Density of $O_M$ in $\mathcal{N}^s_{p,q,r}$}\\
The idea is the same as above.

Observe that the definition implies that there exist integers $n$ and $m$ such that 
\begin{displaymath}
\Big ( \sum_{j \notin \left \{ -n, \dots, 0 \dots m\right\}} 2^{sjr}\vert \vert f_j \vert \vert _{\mathcal{M}^s_{p,q}}^r \Big) ^{1/r} \leq \frac{\varepsilon}{2}.
\end{displaymath}
And as before, this gives us the result that $O_M$ is dense in $\mathcal{N}^s_{p,q,r}$. \\

Another idea to prove the density of $C^{\infty}$ in $N^s_{p,q,r}$ arises from the usual mollification: \\
We have to show that for any given $\varepsilon$ and any given function $f \in N^s_{p,q,r}$ there exists a function $g \in C^{\infty}$ such that 
\begin{displaymath}
\vert \vert f-g \vert N^s_{p,q,r} \vert \vert \leq \varepsilon.
\end{displaymath}
As indicated above, our candidate for $g$ will be a function of the form
\begin{displaymath}
g= \varphi_{\delta} \ast f
\end{displaymath}
where $\varphi_{\delta}$ is a mollifying sequence ( and $\delta$ will be specified later on).\\
First of all, observe that due to Tonelli-Fubini we have $\varphi_{\delta} \ast f^j = (\varphi_{\delta}\ast f)^j$. \\
Now, as above we observe that the fact that $f$ belongs to $N^s_{p,q,r}$ implies that there exists $N_0 \in \mathbb{N}$ such that
\begin{displaymath}
\Big(\sum_{N_0+1}^{\infty} 2^{jsr}\vert\vert f^j \vert M^p_q \vert \vert^r\Big)^{\frac{1}{r}} \leq \tilde{\varepsilon}
\end{displaymath}
which together with \cite{KY}, lemma 1.8, immediately leads to the observation that
\begin{displaymath}
\Big(\sum_{N_0+1}^{\infty} 2^{jsr}\vert\vert (f- f \ast\varphi_{\delta})^j \vert M^p_q \vert \vert^r\Big)^{\frac{1}{r}} \leq \frac{\varepsilon}{2}.
\end{displaymath}
For the remaining contributions we first of all observe that
\begin{displaymath}
\vert f^j - f^j \ast \varphi_{\delta} \vert \leq \vert \vert \nabla f^j \vert \vert_{\infty} \delta \leq C \vert \vert f \vert N^s_{p,q,r} \vert \vert 2^j \delta.
\end{displaymath}
In order to see this, note that $f^j \in N^s_{p,q,1}$ which together with two results from \cite{KY} similar to the estimate (\ref{tao}) and the embedding of Besov-Morrey into Besov spaces (see also \cite{KY}) implies that 
\begin{displaymath}
\vert \vert \nabla f^j \vert \vert_{\infty} \leq C \vert \vert f \vert N^s_{p,q,r} \vert \vert 2^{j}.
\end{displaymath}
In the case $j=0$ observe that 
\begin{eqnarray}
(\partial_{x_i}f)^0&=&\mathcal{F}^{-1}(i\xi_i\hat{f}\phi_0) \nonumber \\
&=& \mathcal{F}^{-1}(i\xi_i\hat{f}\phi_0(\phi_0+\phi_1)) \nonumber \\
&=& f^0 \ast \mathcal{F}^{-1}(i\xi_i(\phi_0+\phi_1)) \nonumber 
\end{eqnarray}
which implies that 
\begin{displaymath}
\vert \vert \partial_{x_i}f^0 \vert M^p_q \vert \vert \leq C \vert \vert f^0\vert M^p_q \vert \vert.
\end{displaymath}
Apart from this observation, the argument is the same as the usual one known in the framework of Lebesgue spaces.\\
Now, we can calculate for any radius $R \in (0,1 \rbrack$ and for any point $x_0 \in \mathbb{R}^n$
\begin{eqnarray}
R^{\frac{n}{p}-\frac{n}{q}}\vert \vert f^j -f^j \ast \varphi_{\delta} \vert \vert_{L^q(B_R(x_0))} 
&=&R^{\frac{n}{p}-\frac{n}{q}}\Big( \int_{B_R(x_0)}\vert f^j -f^j \ast \varphi_{\delta} \vert^q \Big)^{\frac{1}{q}} \nonumber \\
&\leq& C R^{\frac{n}{p}-\frac{n}{q}} \Big( \vert \vert \nabla f^j \vert \vert_{\infty}^q \delta^q R^n\Big)^{\frac{1}{q}} \nonumber \\
&\leq& C R^{\frac{n}{p}-\frac{n}{q}} \Big( \vert \vert f\vert N^s_{p,q,r} \vert \vert^q 2^{jq}\delta^q R^n\Big)^{\frac{1}{q}} \nonumber \\
&=& C R^{\frac{n}{p}} \vert \vert  f\vert N^s_{p,q,r}\vert \vert \delta 2^{j} \nonumber \\
&\leq& C \vert \vert f \vert N^s_{p,q,r}\vert \vert \delta 2^{j} \nonumber 
\end{eqnarray}
from which we conclude that
\begin{eqnarray}
\Big ( \sum_{j=0}^{N_0} 2^{jsr}\vert \vert f^j -f^j \ast \varphi_{\delta} \vert M^p_q \vert \vert^r\Big)^{\frac{1}{r}}
&\leq&  \sum_{j=0}^{N_0} \vert \vert f \vert N^s_{p,q,r} \vert \vert \delta 2^{N_0+N_0sr}\nonumber \\
&\leq& (N_0 +1)  \vert \vert f \vert N^s_{p,q,r} \vert \vert \delta 2^{N_0+N_0sr}\nonumber \\
&\leq& \frac{\varepsilon}{2} \nonumber
\end{eqnarray}
if we choose $\delta$ sufficiently small.\\
This shows that $f \in N^s_{p,q,r}$ can be approximated by compactly supported smooth function - the convolution $f \ast \varphi_{\delta} \ast f$ has compact support.

Now, we assume that $f \in B^s_{\mathcal{M}^p_q,r}$ where $s \geq 0$, $1< q \leq 2$ and $1\leq q \leq p \leq\infty$ has compact support. First of all, we observe that according to lemma \ref{inMorrey} $f \in \mathcal{M}^p_q$ and since it has compact support, $f \in L^q$. From this we deduce that whenever $O\leq j \leq N_0$, $f^j \in B^s_{q,m}$ for all $s \in \mathbb{R}$ and arbitrary $m$ and in particular, $f^j \in L^p$. So for each $j$ there exists a $\delta_j$ such that
\begin{displaymath}
\vert \vert f^j - f^j \ast \varphi_{\delta_j} \vert \vert_q^m \leq \Big(\frac{\varepsilon}{2(N_0+1)}\Big)^m.
\end{displaymath}
If we now choose $\delta$ small enough, then   
\begin{displaymath}
\Big ( \sum_{j=0}^{N_0} 2^{jsr}\vert \vert f^j -f^j \ast \varphi_{\delta} \vert M^p_q \vert \vert^r\Big)^{\frac{1}{r}}
= \Big ( \sum_{j=0}^{N_0} 2^{jsr}\vert \vert (f -f\ast)^j \varphi_{\delta} \vert M^p_q \vert \vert^r\Big)^{\frac{1}{r}}
\leq \frac{\varepsilon}{2}.
\end{displaymath}
The other frequencies are estimated as above.\\
Finally we observe that $f \ast \varphi_\delta$ is not only smooth but also compactly supported since it is a convolution of a compactly supported function with a compactly supported distribution.
\begin{flushright}
$\Box$
\end{flushright}

\begin{Remark}
\textnormal{
A close look at the proof we just gave, shows that in fact
\begin{displaymath}
\cap_{m \geq 0} C^m 
\end{displaymath}
is dense in the above spaces.
}
\end{Remark}

\subsection{Proof of lemma \ref{stab}}

We split the product $fg$ into the three paraproducts $\pi_1(f,g)$, $\pi_2(f,g)$ and $\pi_3(f,g)$ and analyse each of them independently. \\
\begin{itemize}
\item[i)]
We start with $\pi_1(f,g)=\sum_{k=2}^{\infty} \sum_{l=0}^{k-2}f^lg^k$. It is easy to see that a simple adaptation of lemma 3.15 of \cite{Maz2} to our variant of Besov-Morrey, implies that it suffices to show that 
\begin{displaymath}
\Big ( \sum_{k=2}^{\infty} \vert \vert g^k \sum_{l=0}^{k-2}f^l \vert \vert_{\mathcal{M}^n_2}^2 \Big)^{\frac{1}{2}}
\leq C \vert \vert g \vert B^0_{\mathcal{M}^n_2,2} \vert \vert (
\vert \vert f \vert B^1_{\mathcal{M}^n_2,2} \vert \vert + \vert \vert f \vert \vert_{\infty} ).
\end{displaymath}
In fact, we calculate
\begin{eqnarray}
\Big ( \sum_{k=2}^{\infty} \vert \vert g^k \sum_{l=0}^{k-2}f^l \vert \vert^2_{\mathcal{M}^n_2}} \Big)^{\frac{1}{2}&\leq& 
\Big ( \sum_{k=2}^{\infty} \vert \vert g^k (\sup_{s}\vert \sum_{l=0}^{s}f^l\vert ) \vert \vert^2_{\mathcal{M}^n_2}} \Big)^{\frac{1}{2}\nonumber \\
&\leq& \Big ( \sum_{k=2}^{\infty} \vert \vert g^k \vert \vert^2_{\mathcal{M}^n_2} \vert \vert \sup_{s}\vert \sum_{l=0}^{s}f^l\vert \vert \vert_{\infty}^2  \Big)^{\frac{1}{2}}\nonumber \\
&\leq& \vert \vert \sup_{s}\vert \sum_{l=0}^{s}f^l\vert \vert \vert_{\infty} ( \sum_{k=2}^{\infty} \vert \vert g^k \vert \vert^2_{\mathcal{M}^n_2} \Big)^{\frac{1}{2}}\nonumber \\
&\leq& \vert \vert \sup_{s}\vert \sum_{l=0}^{s}f^l\vert \vert \vert_{\infty} \vert \vert g \vert B^0_{\mathcal{M}^n_2,2} \vert \vert \nonumber \\
&\leq& \vert \vert f \vert \vert_{\infty} \vert \vert g \vert B^0_{\mathcal{M}^n_2,2} \vert \vert \nonumber \\
&&\textrm{because of lemma 4.4.2 of \cite{RS}} \nonumber \\
&<&\infty. \nonumber
\end{eqnarray}
\item[ii)]
Next, we study $\pi_2(f,g)=\sum_{k=0}^{\infty} \sum_{l=k-1}^{k+1}f^lg^k$. For our further calculations we fix $l=k$. We will see that what follows will not depend on this choice, so
\begin{displaymath}
\vert \vert \pi_2(f,g) \vert B^0_{\mathcal{M}^n_2,2} \vert \vert \leq C \sup_{s \in \left \{-1,0,1 \right \}} 
\vert \vert \sum_{k=0}^{\infty} f^{k+s}g^k \vert B^0_{\mathcal{M}^n_2,2} \vert \vert.
\end{displaymath}
In fact, we will show a bit more, namely $\pi_2(f,g) \in B^1_{\mathcal{M}^{\frac{n}{2}}_1,1}$. Again a simple adaptation of lemma 3.16 of \cite{Maz2} shows that we only have to estimate $\sum_{k=0}^{\infty} 2^k \vert \vert f^k g^k \vert \vert_{\mathcal{M}^{\frac{n}{2}}_1}$. 
In fact, we have
\begin{eqnarray}
\sum_{k=0}^{\infty} 2^k \vert \vert f^k g^k \vert \vert_{\mathcal{M}^{\frac{n}{2}}_1} &\leq& \sum_{k=0}^{\infty} 2^k \vert \vert f^k\vert \vert_{\mathcal{M}^n_2} \vert \vert g^k \vert \vert_{\mathcal{M}^n_2} \nonumber \\
&\leq& \Big( \sum_{k=0}^{\infty} 2^{2k} \vert \vert g^k\vert \vert^2_{\mathcal{M}^n_2} \Big)^{\frac{1}{2}} \Big(\sum_{k=0}^{\infty} \vert \vert f^k \vert \vert^2_{\mathcal{M}^n_2}\Big)^{\frac{1}{2}} \nonumber \\
&\leq& \vert \vert g \vert B^1_{\mathcal{M}^n_2,2} \vert \vert \; \vert \vert f \vert B^0_{\mathcal{M}^n_2,2} \vert \vert \nonumber \\
&<&\infty. \nonumber
\end{eqnarray}
Once we have this, it implies together with the embedding of Besov-Morrey spaces into Besov spaces (see \cite{KY}) - adapted to our variant of Besov-Morrey spaces - and the fact that $l^1 \subset l^2$ immediately that $\sum_{k=0}^{\infty} f^k g^k \in B^0_{\mathcal{M}^n_2,2}$. And finally we get that $\pi_2(f,g) \in B^0_{\mathcal{M}^n_2,2}$.
\item[iii)]
The remaining addend is $\pi_3(f,g)$. Again, as in i) it is enough to show that we can estimate $\Big ( \sum_{l=2}^{\infty} \vert \vert f^l \sum_{k=0}^{l-2}g^k \vert \vert_{\mathcal{M}^n_2}^2 \Big)^{\frac{1}{2}}$ in the desired manner. In fact we observe that the following inequalities hold:
\begin{eqnarray}
\Big ( \sum_{l=2}^{\infty} \vert \vert f^l \sum_{k=0}^{l-2}g^k \vert \vert_{\mathcal{M}^n_2}^2 \Big)^{\frac{1}{2}} &\leq&
\sum_{l=2}^{\infty} \vert \vert f^l \sum_{k=0}^{l-2}g^k \vert \vert_{\mathcal{M}^n_2} \nonumber \\
&\leq& \sum_{l=2}^{\infty} \vert \vert f^l \vert \vert_{\mathcal{M}^n_2} \vert \vert\sum_{k=0}^{l-2}g^k \vert \vert_{\infty} \nonumber \\
&=& \sum_{l=2}^{\infty} 2^l\vert \vert f^l \vert \vert_{\mathcal{M}^n_2} 2^{-l}\vert \vert\sum_{k=0}^{l-2}g^k \vert \vert_{\infty} \nonumber \\
&\leq& \Big (\sum_{l=0}^{\infty} 2^{2l} \vert \vert f^l \vert \vert^2_{\mathcal{M}^n_2}\Big )^{\frac{1}{2}} \Big ( \sum_{l=0}^{\infty} 2^{-2l} \vert \vert\sum_{k=0}^{l-2}g^k \vert \vert^2_{\infty} \Big)^{\frac{1}{2}} \nonumber \\
&\leq& C\Big (\sum_{l=0}^{\infty} 2^{2l} \vert \vert f^l \vert \vert^2_{\mathcal{M}^n_2}\Big )^{\frac{1}{2}} \Big ( \sum_{l=0}^{\infty} 2^{-2l} \vert \vert\sum_{k=0}^{l}g^k \vert \vert^2_{\infty} \Big)^{\frac{1}{2}} \nonumber \\
&\leq& C\vert \vert f \vert B^1_{\mathcal{M}^n_2,2}\vert \vert \Big ( \sum_{l=0}^{\infty} 2^{-2l} \vert \vert\sum_{k=0}^{l}g^k \vert \vert^2_{\infty} \Big)^{\frac{1}{2}} \nonumber \\
&\leq& C\vert \vert f \vert B^1_{\mathcal{M}^n_2,2} \vert \vert \;  \vert \vert g \vert B^{-1}_{\infty,2}\vert \vert \nonumber \\
&&\textrm{according to lemma 4.4.2 of \cite{RS}} \nonumber \\
&\leq& C\vert \vert f \vert B^1_{\mathcal{M}^n_2,2} \vert \vert \; \vert \vert g \vert N^0_{n,2,2}\vert \vert \nonumber \\
&&\textrm{due to the embedding result for Besov-Morrey spaces (\cite{KY})} \nonumber \\
&\leq& C\vert \vert f \vert B^1_{\mathcal{M}^n_2,2} \vert \vert \;  \vert \vert g \vert B^0_{\mathcal{M}^n_2,2}\vert \vert \nonumber \\
&<&\infty. \nonumber
\end{eqnarray}
\end{itemize}
If we put together all our results from i) to iii) we see that we have the estimate
\begin{displaymath}
\vert \vert gf \vert B^0_{\mathcal{M}^n_2,2} \vert \vert \leq C \vert \vert g \vert B^0_{\mathcal{M}^n_2,2} \vert \vert (
\vert \vert f \vert B^1_{\mathcal{M}^n_2,2} \vert \vert + \vert \vert f \vert \vert_{\infty} )
\end{displaymath}
as claimed.
\begin{flushright}
$\Box$
\end{flushright}


\end{document}